\documentstyle[twoside,12pt]{article}

\input amssym.def
\input amssym

\newcount\probno
\newcount\sectorno
\def\sector#1{\probno=1\advance\sectorno by 1
 \newpage\begin{center}{\bf #1}\end{center}}
\def\prob{\par\medskip\noindent
 {\kern -.5in \parbox{.45in}{\the\sectorno .\the\probno}}\advance\probno by 1}
\newtheorem{lemma}{Lemma}

\newwrite\dexno
\immediate\openout\dexno=dexdef
\def\dex#1{{\it #1}\immediate\write\dexno{?par?noindent #1 ?dotfill
   \the\sectorno.\the\probno -p.\thepage}}
\def\sdex#1{\immediate\write\dexno{?par?noindent #1 ?dotfill
   \the\sectorno.\the\probno -p.\thepage}}

\def\ax #1.{\par\medskip\noindent{\bf #1:}}
\def\blank{\;\tilde{ }\;}   \def\cc{{\frak c}}  \def\conj{\mathop{\wedge}}
\def\disj{\mathop{\vee}}    \def\dom{{\rm dom}}  \def\dotminus{\dot{-}}
\def\elemsub{\preceq}       \def\elemsup{\succeq} \def\eq{\approx}
\def\fama{{\cal A}}         \def\famc{{\cal C}}   \def\famf{{\cal F}}
\def\HFin{(HF,\in)}         \def\implies{\rightarrow} \def\isom{\simeq}
\def\moda{{\goth A}}        \def\modb{{\goth B}}  \def\modc{{\goth C}}
\def\modn{{\goth N}}        \def\modq{{\goth Q}}  \def\modr{{\goth R}}
           \def\nor{\oplus}      \def\ord{{\bf ORD }}
\def\partial{{\rm partial}{\;\;}}
\def\pp{{\cal P}} \def\prf{{\rm Prf}}
\def\proof{\par\noindent proof:\par}              \def\proves{\vdash}
\def\qed{\nopagebreak\par\noindent\nopagebreak$\Box$\par}
\def\qq{{\Bbb Q}}          \def\res{\upharpoonright}
\def\rmand{{\mbox{ and }}} \def\rmiff{{\mbox{ iff }}}
\def\rmor{{\mbox{ or }}}   \def\rr{{\Bbb R}}       \def\ss{{\cal S}}
\def\substruc{\subseteq}   \def\ul{\ulcorner}      \def\ur{\urcorner}
\def\zfc{{\rm ZFC}}        \def\zz{{\Bbb Z}}
\def\com#1{}

\begin{document}

\begin{flushright}
  Arnold W. Miller\\
  Department of Mathematics\\
  University of Wisconsin\\
  480 Lincoln Drive\\
  Madison, WI 53706\\
  miller@math.wisc.edu\\
  Fall 95\\
\end{flushright}

\bigskip

\begin{center}
  Introduction to Mathematical Logic
\end{center}

I have used these questions or some variations four times to
teach a beginning graduate course in Mathematical Logic.  I
want to thank the many students who hopefully
had some fun doing them,  especially,  Michael Benedikt,
Tom Linton, Hans Mathew, Karl Peters, Mark Uscilka, Joan Hart,
Stephen Mellendorf, Ganesan Ramalingam, Steven Schwalm, Garth
Dickie, Garry Schumacher, Krung Sinapiromsaran, Stephen Young,
Brent Hetherwick, Maciej Smuga-Otto, and Stephen Tanner.

\begin{center}
  Instructions
\end{center}

Do not read logic books during this semester, it is self-defeating.
You will learn proofs you have figured out yourself and the more
you have to discover yourself the better you will learn them.
You will probably not learn much from your fellow student's presentations
(although the one doing the presenting does).  And you shouldn't!
Those that have solved the problem should be sure that the presented
solution is correct.  If it doesn't look right it probably isn't.
Don't leave this up to me, if I am the only one who objects I will
stop doing it.  For those that haven't solved the problem, you
should regard the presented solution as a hint and go and write
up for yourself a complete and correct solution.  Also you might want
to present it to one of your fellow students outside the classroom, if
you can get one to listen to you.

\newpage

\begin{center}
   The Moore Method
\end{center}

From P.R. Halmos\footnote{The teaching of problem solving,
Amer. Math. Monthly, (1975)82, 466-470.}:

``What then is the secret--what is the best way to learn to solve
problems?  The answer is implied by the sentence I started with: solve
problems. The method I advocate is sometimes known as the `Moore
method', because R.L. Moore developed and used it at the University
of Texas.  It is a method of teaching, a method
of creating the problem-solving attitude in a student, that
is a mixture of what Socrates taught
us and the fiercely competitive spirit of the Olympic games.''

\bigskip

From F.Burton Jones\footnote{The Moore method,
Amer. Math. Monthly, (1977)84, 273-278.}:

``What Moore did: $\ldots$
After stating the axioms and giving motivating examples to illustrate
their meaning he would then state definitions and theorems.  He simply
read them from his book as the students copied them down.
 He would then instruct the class to find proofs of their own and
 also to construct examples to show that the hypotheses of the theorems
 could not be weakened, omitted, or partially omitted.

$\ldots$

``When a student stated that he could prove Theorem x, he was asked to
go to the blackboard and present the proof.  Then the other students,
especially those who hadn't been able to discover a proof, would
make sure that the proof presented was correct and convincing.
Moore sternly prevented heckling.  This was seldom necessary because
the whole atmosphere was one of a serious community effort to understand
the argument.''

\bigskip

From D.Taylor\footnote{Creative teaching: heritage of R.L.Moore, University of
Houston, 1972, QA 29 M6 T7, p149.}:

``Criteria which characterize the Moore method of teaching include:

\noindent (1) The fundamental purpose: that of causing a student to
develop his power at rational thought.

\noindent (2) Collecting the students in classes with common mathematical
knowledge,
striking from membership of a class any student whose knowledge is too
advanced over others in the class.

\noindent (3) Causing students to perform research at their level by
confronting
the class with impartially posed questions and conjectures which are at
the limits of their capability.

\noindent (4) Allowing no collective effort on the part of the students
inside
or outside of class, and allowing the use of no source material.

\noindent (5) Calling on students for presentation of their efforts at
settling
questions raised, allowing a feeling of ``ownership'' of a theorem
to develop.

\noindent (6) Fostering competition between students over the settling
of questions raised.

\noindent (7) Developing skills of critical analysis among the class by
burdening
students therein with the assignment of ``refereeing'' an argument
presented.

\noindent (8) Pacing the class to best develop the talent among its
membership.

\noindent (9) Burdening the instructor with the obligation to not assist, yet
respond to incorrect statements, or discussions arising from incorrect
statements, with immediate examples or logically sound propositions to
make clear the objection or understanding.''

\bigskip

Taylor's (2) and (4) are a little too extreme for me.
It is OK to collaborate
with your fellow students as long as you give them credit.  In fact, it
is a good idea to try out your argument first by presenting it to
fellow student.   Avoid reading logic if you can, at least this semester,
but if you do give a reference.

For more readings on the Moore method see:

Paul R. Halmos, What is Teaching?,
Amer. Math. Monthly, 101 (1994), 848-854.

Donald R. Chalice, How to teach a class by the modified Moore method,
Amer. Math. Monthly, 102 (1995), 317-321.

\bigskip

\begin{center}
Quote From P.R. Halmos:
\end{center}

``A famous dictum of P\'{o}lya's about problem solving is that
if you can't solve a problem, then there is an easier problem that
you can't solve--find it!''

\sector{Propositional Logic and the Compactness Theorem}

The \dex{syntax} (grammar) of propositional logic is the following.
The \dex{logical symbols} are $ \conj ,\disj , \neg , \implies ,$ and $ \iff$.
The \dex{nonlogical symbols} consist of an arbitrary nonempty set
$\pp$ that we assume
is disjoint from the set of logical symbols to avoid confusion.  The
set ${\pp}$ is referred to as the set of \dex{atomic sentences} or
as the set of \dex{propositional letters}.
For example, $\{P,Q,R\}$, $\{P_0,P_1,P_2,\ldots\}$, or $\{S_r:r\in\rr\}$.
The set of
\dex{propositional sentences} ${\ss}$ is the smallest set of finite strings
of symbols such that
$\pp\subseteq \ss$, and if $\theta \in \ss$ and $\psi \in \ss$, then
$ \neg \theta \in \ss, (\theta \conj \psi) \in \ss,(\theta \disj \psi) \in \ss,
(\theta \implies \psi) \in \ss,$ and $(\theta \iff \psi) \in \ss$.

The \dex{semantics} (meaning) of propositional logic consists of truth
 evaluations.
A \dex{truth evaluation} is a function $e:{\ss}\to\{ T,F \}$,
that is consistent
with the following truth tables:
$$
\begin{array}{ccccccc}
\theta & \psi &\neg\theta &(\theta\conj\psi) &(\theta\disj\psi)
&(\theta\implies\psi) &(\theta\iff\psi) \\
T&T&F&T&T&T&T\\
T&F&F&F&T&F&F\\
F&T&T&F&T&T&F\\
F&F&T&F&F&T&T
\end{array}
$$
For example if $e(\theta)=T$ and $e(\psi)=F$, then
$e(\theta \implies \psi)=F$. Also $e(\neg\theta)=T$ iff $e(\theta)=F$.
For example, if $\pp=\{P_x:x\in\rr\}$ and we define
$e(P_x)=T$ if $x$ is a rational and  $e(P_x)=T$ if $x$ is a irrational,
then $e((P_2\conj \neg P_{\sqrt{2}}))=T$.  However if we
define $e^\prime(P_x)=T$ iff $x$ is an algebraic number, then
$e^{\prime}((P_2\conj \neg P_{\sqrt{2}}))=F$.

A sentence
$\theta$ is called a \dex{validity} iff for every truth evaluation $e$,
$e(\theta)=T$.  A sentence
$\theta$ is called a \dex{contradiction} iff for every truth evaluation $e$,
$e(\theta)=F$.

We say that two sentences $\theta$ and $\psi$ are \dex{logically equivalent}
iff for every truth evaluation $e$, $e(\theta)=e(\psi)$.
A set of logical symbols is \dex{adequate} for propositional logic iff
every propositional sentence is logically equivalent to one whose only logical
symbols are from the given set.

\bigskip

\prob Define ${\ss}_0={\pp}$ the atomic sentences and define
$${\ss}_{n+1}={\ss}_n\cup\{\neg\theta:\theta\in {\ss}_n\}\cup
\{(\theta\#\psi):\theta,\psi\in {\ss}_n,\#\in
\{\conj,\disj,\implies,\iff\}\}$$
Prove that ${\ss}={\ss}_0\cup {\ss}_1\cup {\ss}_2\cup\cdots$.

\prob Prove that for any function $f:{\pp}\to \{T,F\}$ there exists a
unique truth evaluation $e:\ss\to \{T,F\}$ such that  $f=e\res {\pp}$.
The symbol $e\res {\pp}$
stands for the restriction of the function $e$ to $\pp$.
\sdex{$e?res {?pp}$}

\prob Let $\theta$ and $\psi$ be two propositional sentences.  Show
that $\theta$ and $\psi$ are logically equivalent iff $(\theta\iff\psi)$
is a validity.

\prob Suppose $\theta$ is a propositional validity, $P$ and $Q$ are
two of the propositional letters occurring in $\theta$, and $\psi$
is the sentence obtained by replacing each occurrence of $P$ in
$\theta$ by $Q$.   Prove that $\psi$ is a validity.

\prob Can you define $ \disj$  using only $\implies$?  Can you define $\conj$
 using only $\implies$?

\com{ $(p\implies q)\implies p$ is equivalent to $p\disj q$. $conj$ cannot
be defined (if there are at least two atomic sentences- it is definable if
there is only one), one proof is that every sentence using only implies is
made true by at least 50\% of the truth evaluations.}

\prob Show that  $\{ \disj, \neg \}$ is an adequate set for propositional
logic.

\prob The definition of the logical connective \dex{nor} ( $\nor$ )
is given by
the following truth table:
$$
\begin{array}{ccc}
 \theta & \psi & (\theta \nor \psi) \\
           T&T&F  \\
           T&F&F  \\
           F&T&F  \\
           F&F&T
\end{array}
$$
Show that $\{ \nor \}$ is an adequate set for propositional logic.

\prob (Sheffer) Find another binary connective that is adequate all
by itself.

\prob Show that   $\{\neg\}$   is not adequate.

\prob Show that $\{ \disj \}$ is not adequate.

\prob How many binary logical connectives are there?  We assume two
connectives are the same if they have the same truth table.

\prob Show that there are exactly two binary logical connectives that are
 adequate
all by themselves. Two logical connectives are the same iff they have the same
truth table.

\prob Suppose $\pp=\{P_1,P_2,\ldots,P_n\}$. How many propositional
sentences (up to logical equivalence) are there in this language?

\prob Show that every propositional sentence is equivalent to a sentence in
\dex{disjunctive normal form}, i.e. a disjunction of conjunctions
of atomic or the negation of atomic sentences:
$$\disj_{i=1}^m \bigl(\conj_{j=1}^{k_i}\theta_{ij}\bigr)$$
where each $\theta{ij}$ is atomic or $\neg$atomic.
The expression $\disj_{i=1}^n \psi_i$  abbreviates
$(\psi_1\disj(\psi_2\disj(\cdots\disj (\psi_{n-1}\disj\psi_{n})))\cdots)$.

\bigskip
In the following definitions and problems $\Sigma$ is a set of propositional
sentences in some fixed language and all sentences are assumed to be in this
same fixed language.
$\Sigma$ is \dex{realizable} iff there exists a truth evaluation $e$ such that
for all $\theta \in \Sigma$,  $e(\theta)=T.$
$\Sigma$ is \dex{finitely realizable}
iff every finite subset of $\Sigma$ is realizable.
$\Sigma$ is \dex{complete} iff for every sentence $\theta$ in the language of
$\Sigma$ either  $\theta$ is in $\Sigma$ or $\neg\theta$ is in $\Sigma$.

\medskip

\prob Show that if $\Sigma$ is finitely realizable and $\theta$ is
any sentence
 then either
$\Sigma \cup \{ \theta \}$ is finitely
realizable or $\Sigma \cup \{ \neg\theta \}$
 is finitely realizable.

\prob Show that if $\Sigma$ is finitely realizable and  $(\theta \disj \psi)$
 is in $\Sigma$, then either
$\Sigma \cup \{ \theta\}$ is finitely realizable or
$\Sigma \cup \{ \psi \}$ is finitely realizable.

\prob Show that if $\Sigma$ is finitely realizable
and complete and if $\theta$ and
$(\theta \implies \psi)$
are both in $\Sigma$, then $\psi$ is in $\Sigma$.

\prob Show that if $\Sigma$ is finitely realizable
and complete, then $\Sigma$ is
realizable.

\prob Suppose that the set of all sentences in our language is countable,
e.g.,
$S =\{ \theta_n : n = 0,1,2,\ldots \}$. Show that if $\Sigma$ is finitely
 realizable, then there
exists a complete finitely
realizable $\Sigma^\prime $ with $\Sigma\subseteq \Sigma^\prime$.

\prob ({\bf Compactness theorem for propositional logic})
Show that every finitely realizable $\Sigma$ is realizable.  You may assume
there are only countably many sentences in the language.

\bigskip
 A family of sets $\famc$ is a \dex{chain} iff for any $X,Y$ in $\famc$
 either
 $X \subseteq Y$ or $Y \subseteq X$.
The union of the family $\fama$ is
$$\bigcup \fama = \{ b : \exists c \in \fama ,b \in c \}.$$
$M$ is a \dex{maximal} member of a family $\fama$ iff $M \in \fama$ and
 for every $B$ if
$B \in \fama$ and $M \subseteq B$, then $M=B$.
A family of sets $\fama$ is closed under the unions of chains iff
for every subfamily, $\famc$,
of $\fama$ which is a chain the union of the chain, $\bigcup\famc$,
is also a member of $\fama$.

\par\medskip
 {\bf Maximality Principle:}  Every family of sets closed
 under the unions of
chains has a maximal member.
\par\medskip

\prob Show that the family of finitely realizable $\Sigma$ is closed under
unions of chains.

\prob Show how to prove the compactness theorem without the assumption
 that there are only countably many sentences. (You may use the Maximality
 Principle.)

\prob Suppose $\Sigma$ is a set of sentences and $\theta$ is some sentence
such that for every
truth evaluation $e$  if $e$ makes all sentences in $\Sigma$ true, then
$e$ makes
$\theta$ true.
Show that for some
 finite $\{ \psi_1, \psi_2, \psi_3,\ldots \psi_n \}\subseteq \Sigma$
the sentence
$$ ( \psi_1 \conj \psi_2 \conj \psi_3 \conj\cdots \conj \psi_n )
\implies \theta$$
is a validity.

\par\bigskip
\noindent A \dex{binary relation} $R$ on a set $A$ is a subset of $A \times A$.
Often
we write $xRy$ instead of $\langle x,y\rangle\in R$.
 A binary relation $\leq$ on a set $A$ is a \dex{partial order} iff
\par   a. (reflexive) $\forall a \in A \;  a\leq a$;
\par   b. (transitive) $ \forall a,b,c \in A \;[ ( a \leq b  \conj  b \leq c)
\implies  a \leq c]$; and
\par   c. (antisymmetric) $\forall a,b\in A\;[ (a \leq b  \conj  b \leq a)
\implies  a=b]$.

\noindent Given a partial order $\leq$ we define the \dex{strict order}
 $<$ by
$$x < y \iff (x \leq y \conj  x \not= y)$$

\noindent A binary relation $\leq$ on a set A is a \dex{linear order}
 iff $\leq$
is a partial order and
\par   d. (total) $\forall a,b \in A (a \leq b \disj b \leq a)$.

  A binary relation $R$ on a set $A$ extends a binary relation
$S$ on $A$ iff $S \subseteq R$.

\medskip

\prob Show that for every finite set $A$ and partial order $\leq$ on $A$
 there exists a
linear order $\leq ^*$ on $A$ extending $\leq$.

\prob Let $A$ be any set and let our set of atomic sentences $\pp$ be:
        $$\pp= \{ P_{ab} : a,b \in A \}$$
For any truth evaluation $e$ define  $\leq_e$ to be the binary relation
on $A$ defined by
       $$ a \leq_e b   \mbox{ iff }  e(P_{ab})=T. $$
Construct a set of sentences $\Sigma$ such that for every truth
evaluation $e$,

\centerline{$e$ makes $\Sigma$ true iff $\leq_e$ is a linear order on $A$.}

\prob  Without assuming the set $A$ is finite prove
 for every partial order $\leq$ on $A$
 there exists a linear order $\leq ^*$ on $A$ extending $\leq$.
\com{This one can be proved directly from the maximality principal, the
others I don't know about.}

\medskip
 In the next problems $n$ is an arbitrary positive integer.
\par\medskip

\prob If $X \subseteq  A$ and $R$ is a binary relation on $A$ then the
restriction of $R$ to $X$ is the binary relation
$S= R\cap (X \times X)$.
For a partial order $\leq$ on $A$, a set $B \subseteq A$ is an $\leq$-chain
iff the
restriction of $\leq$ to $B$ is a linear order.
 Show that given a  partial order $\leq$ on $A$:

the set $A$  is the union of less than $n$  $\leq$-chains
  iff
every finite subset of $A$ is the union of less than $n$  $\leq$-chains.

\prob A partial order $\leq$ on a set $A$ has \dex{dimension}
 less than $n+1$ iff
 there exists $n$ linear
orders $\{ \leq_1, \leq_2, \leq_3,\ldots,\leq_{n} \}$ on $A$ (not
necessarily distinct) such that:
       $$ \forall x,y \in A\;  [  x\leq y   \mbox{ iff }
          ( x\leq_i y \mbox{ for } i=1,2,\ldots,n) ].$$
Show that a partial order $\leq$ on a set $A$ has dimension less than $n+1$
iff for every
finite $X$ included in $A$ the restriction of $\leq$ to
 $X$ has dimension less than $n+1$.

\prob
A binary relation $E$ (called the edges) on a set $V$ (called the vertices)
is a \dex{graph} iff
\par    a. (irreflexive) $\forall x \in V    \neg xEx $; and
\par    b. (symmetric) $\forall x,y \in V \,  (xEy \implies yEx)$.
\par\noindent  We say $x$ and $y$ are adjacent iff $xEy$.
$(V^\prime ,E^\prime )$ is a subgraph of $(V,E)$ iff
$V^\prime\subseteq V$ and $E^\prime$ is the restriction of $E$ to $V^\prime$.
For a graph $(V,E)$ an $n$ coloring is a map
$c : V \to \{ 1,2,\ldots,n \}$ satisfying
$\forall x,y \in V  (xEy \implies c(x) \not= c(y))$,
 i.e. adjacent vertices have different colors.
A graph $(V,E)$ has \dex{chromatic number} $\leq n$ iff there is a $n$
coloring on its vertices.
 Show that a graph has chromatic number $\leq n$ iff
  every finite subgraph of it has
chromatic number $\leq n$.

\prob A triangle in a graph $(V,E)$ is a set $\Delta=\{a,b,c\} \subseteq V$
such that $aEb$, $bEc$, and $cEa$.  Suppose that every finite subset of
$V$ can be partitioned into $n$ or fewer sets none of which contain a
triangle.
Show that $V$ is the union of $n$ sets none of
which contain a triangle.

\prob (Henkin) A \dex{transversal} for a family of sets $\famf$ is a
one-to-one choice function.
That is a one-to-one function $f$ with domain $\famf$ and for every
$x\in \famf \,\, f(x)\in x.$ Suppose that $\famf$ is a family of
finite sets such that for
every finite $\famf^\prime  \subseteq \famf, \famf^\prime $ has a
transversal.
Show that $\famf$ has a
transversal.  Is this result true if $\famf$ contains infinite sets?

\prob Let $\famf$ be a family of subsets of a set $X$.  We say
that $\famc\subseteq \famf$
is an \dex{exact cover}
 of $Y \subseteq X$ iff every element of $Y$ is in a unique
element of $\famc$.  Suppose that every element of $X$ is in at most
finitely many
elements of $\famf$.  Show that there exists an exact cover
$\famc\subseteq \famf$ of $X$
iff for every finite $Y\subseteq X$ there exists
$\famc\subseteq \famf$ an exact
cover
of $Y$. Is it necessary that every element of $X$ is in at most finitely many
elements of $\famf$?

\com{Use $P_a$ to say ``$a\in \famc$'', so $\Sigma$ contains
$\disj_{x\in a} P_a$ for each $x\in X$ and $\neg(P_a\conj P_b)$ if
$a$ and $b$ are not disjoint. To verify finite sat- given Y finite get
$Z\supseteq Y$ finite which can witness the nondisjointedness of the
elements of $\famf$ which hit $Y$, then for exact cover $Z$ the
elements which
hit $Y$ are a disjoint exact cover of $Y$.}

\prob  If $\famf$ is a family of subsets of $X$ and $Y\subseteq X$
then we say $Y$
\dex{splits} $\famf$ iff for any $Z\in \famf,$
$ Z \cap Y $ and $Z \setminus Y$
 are both nonempty. Prove that if $\famf$ is a family of finite
subsets of $X$
then
$\famf$ is split by some $Y\subseteq X$ iff every finite
$\famf^\prime \subseteq \famf$
is split by some $Y\subseteq X$.  What if $\famf$ is allowed to
have infinite sets
in it?

\prob Given a set of students and set of classes, suppose each student
wants one of a finite set of classes, and each class has a finite
enrollment limit.  Show that if each finite set of students can
be accommodated, they all can be accommodated.

\prob Show that the compactness theorem of propositional logic is
equivalent to the statement that for any set $I$, the space $2^I$, with
the usual Tychonov product topology is compact, where $2=\{0,1\}$ has
the discrete topology. (You should skip this problem if you do not know what
a topology is.)

\sector{The Axioms of Set Theory}

   Here are some.  The whole system is known as \dex{ZF} for Zermelo-Fraenkel
   set theory.  When the axiom of choice is included it is denoted
   \dex{ZFC}.  It was originally developed by Zermelo to
   make precise what he meant when he said that the
   well-ordering principal follows from the
   axiom of choice.   Latter Fraenkel added the axiom of replacement.
   Another interesting system is GBN which is G\"{o}del-Bernays-von Neumann
   set theory.

\medskip

\ax Empty Set. $$ \exists x\, \forall y ( y \notin x )$$
The empty set is usually written $\emptyset$. \sdex{$?emptyset$}

\ax Extensionality.
 $$ \forall x \forall y ( x = y \iff \forall z ( z \in x \iff z\in y)) $$
Hence there is only one empty set.

\ax Pairing. $$ \forall x \forall y \exists z \forall u
( u \in z \iff u = x \disj u = y ) $$
We usually write $z = \{ x,y \}$. \sdex{$z = ?{ x,y ?}$}

\ax Union. $$ \forall x\, \exists y\, ( \forall z ( z \in y \iff
(\exists u\, u \in x \conj z \in u ))$$
We usually write $y = \cup x$. \sdex{$y = \cup x $}  $A \cup B$
abbreviates  $ \cup \{ A,B \}$.  \sdex{$A \cup B$}
$z \subseteq x $  is an abbreviation
for $ \forall u \, (u\in z \implies u\in x)$. \sdex{$z ?subseteq x $}

\ax Power Set. $$ \forall x \, \exists y \, \forall z
( z \in y \iff z \subseteq x ) $$
We usually write $y = P(x)$. \sdex{$y = P(x)$}
 For any set x , $x + 1 = x \cup \{x\}$. \sdex{$x + 1 = x ?cup ?{x?}$}

\ax Infinity. $$\exists y \, (\emptyset\in y \conj \forall x
( x \in y \implies x + 1 \in y )) $$
The smallest such y is denoted $\omega$, \sdex{$?omega$} so
$\omega=\{0,1,2,\ldots\}$.

\ax Comprehension Scheme.
$$ \forall z \exists y \forall x [ x \in y \iff
                                 ( x \in z \conj \theta(x) ) ]$$
The comprehension axiom is being invoked when we say given $z$ let
$$y = \{ x \in z : \theta(x) \}.$$ \sdex{$y = ?{ x ?in z : ?theta(x) ?}$}
The formula $\theta$ may refer to $z$ and
to other sets, but not to $y$.  In general given a formula $\theta(x)$
the family $\{x:\theta(x)\}$ is referred to as a \dex{class}, it may not
be a set.    For example, the class of all sets is
$${\bf V}=\{x:x=x\}.$$ \sdex{$V$}
Classes that are not sets are referred to as
\dex{proper classes}.
Every set $a$ is a class, since the formula
``$x\in a$'' defines it.
The comprehension axioms say that the intersection
of a class and a set is a set.  We use {\bf boldface} characters to
refer to classes.

\bigskip

\prob Define  $X\cap Y$,
$X\setminus Y$, and
$\bigcap X$ and show they exist. \sdex{$X?cap Y$}
\sdex{$X?setminus Y$}  \sdex{$?bigcap X$}

\prob The \dex{ordered pair} is defined by
\sdex{$?langle x,y?rangle$}
$${\langle x,y\rangle =\{\{x\},\{x,y\}\}}.$$
Show it exists.
Show the $\langle x,y\rangle = \langle u,v\rangle$ iff $x = u$ and $y = v$.

\prob The \dex{cartesian product} is defined by
\sdex{$X ?times Y$}
$$ X \times Y =\{ \langle x,y\rangle : x \in X \mbox{ and } y \in Y\}.$$
 Show it exists.

\prob A function is identified with its graph. For any sets $X$ and $Y$ we let
$Y^X$ be the set of all functions with domain $X$ and
range $Y$.  \sdex{$Y^X$} Show this
set exists.

\prob Given a function $f:A\mapsto B$ and set $C\subseteq A$ the
\dex{restriction} of $f$ to $C$, written $f\res C$
is the function
with domain $C$
and equal to $f$ everywhere in $C$.  \sdex{$f?res C$} Show that it exists.
$f^{\prime\prime}C$   is the set
of all elements of $B$ that are in the image of  $C$.
\sdex{$f^{\prime\prime}C$} Show that it exists.

\prob Prove $\omega$ exists (i.e.
that there does exist a smallest such $y$).
Prove for any formula $\theta(x)$ if $\theta(0)$ and
$\forall x\in\omega (\theta(x)\implies \theta(x+1))$, then
$\forall x\in\omega\;\theta(x)$.

\prob Suppose $G:Z\to Z$.
Show that for any $x\in Z$ there exists a unique $f:\omega\to Z$
such that $f(0)=x$ and for all $n\in\omega\;\; f(n+1)=G(f(n))$.

\prob Let $(V,E)$ be a graph.  Informally, two vertices in any graph
are \dex{connected} iff (either they are the same or)  there is a
finite path using the edges of the graph connecting one to the other.
Use the preceding problem to formally define and prove that
the relation $x$ is connected to $y$, written $x\sim y$,
exists and is an equivalence relation on $V$.
Equivalence classes of $\sim$ are called the
\dex{components} of the graph.

\prob Let $A$ and $B$ be disjoint sets and $f:A\to B$ and $g:B\to A$
be one-to-one functions.  Consider the \dex{bipartite graph} $V$ which
has vertices $A\cup B$ and edges given by the union of the graphs of
$f$ and $g$, i.e., there is an edge between $a\in A$ and $b\in B$ iff
either $f(a)=b$ or $g(b)=a$.
Describe what the finite components of $V$ must look like as a
subgraph of $V$.
Describe the infinite components of $V$.

\prob Define $| X | = | Y |$  iff
there is a one-to-one onto map from $X$ to $Y$.  \sdex{$| X | = | Y |$}
We say $X$ and $Y$ have the same cardinality.
Define   ${| X | \leq | Y |}$ iff
there is a one-to-one map
from $X$ to $Y$.  \sdex{$| X | ?leq | Y |$}
Define ${| X | < |Y |}$ iff $ | X | \leq | Y |$ and $ | X | \not= | Y |$.
\sdex{$| X | < |Y |$}
(Cantor-Shr\"oder-Bernstein)
Show that if $| A | \leq | B |$ and $| B | \leq | A |$ then
$| A | = | B |$.

\bigskip

\prob Show that $|A|\leq |A|$.
Show that if $| A | \leq | B |$ and $| B | \leq | C |$ then
$| A | \leq | C |$.

\prob Show that
$|P(X)| = |\{ 0,1\}^X|$.

\prob (Cantor) Show that  $| X | < | P(X) |$.

\prob Show that the class of sets, ${\bf V}$, is not a set.

\prob Show that $|A \times ( B \times C )| = | ( A \times B ) \times C |$.

\prob Show that  $| A^{ B \times C}| =  | {( A^B)}^C |$.

\prob Show that if there is a function $f:A\mapsto B$ that is onto,
then $|B|\leq |A|$.
\footnote{This requires the axiom of choice. It is open if it is equivalent
to AC.}

\bigskip

A set is finite iff it can be put into one-to-one correspondence with
an element of $\omega$
A set is \dex{countable} iff it is either finite or of the same cardinality
as $\omega$. A set is \dex{uncountable} iff it is not countable.
${{\rr}}$   is the set of real numbers and we use
$\cc = |\rr|$
to denote its cardinality which is also called
the cardinality of the continuum. \sdex{${?rr}$} \sdex{$?cc$}
Below, you may use whatever set theoretic
definitions of the integers, rationals and real numbers that you know.
For example, you may regard the reals as either Dedekind cuts in
the rationals, equivalences classes of Cauchy sequences of rationals,
infinitely long decimals, or ?points on a line.

\bigskip

\prob Show that the set of integers $\zz$
is countable.\sdex{${?zz}$}

\prob Show that the set of odd positive integers is countable.

\prob Show that the set of points in the plane with integer coordinates is
   countable.

\prob Show that the countable union of countable sets is countable.
\footnote{Do you think you needed to use the Axiom of Choice?}

\prob
For any set $X$ let $[X]^{<\omega}$
be the finite subsets of $X$.  \sdex{$[X]^{<?omega}$}
Show that the set of finite subsets of $\omega$, which is
written $[\omega]^{<\omega}$,
is countable.  \sdex{$[?omega]^{<?omega}$}

\prob Show that if there are only countably many atomic sentences then the
set of all propositional sentences is countable.

\prob Show that the set of rationals $\qq$
is countable. \sdex{${?qq}$}

\com{
A clever mapping of the positive rationals to the positive
integers is ${p_1^{k_1}\cdots p_n^{k_n}}\over {q_1^{l_1}\cdots q_m^{l_m}}$
is mapped to ${p_1^{2k_1}\cdots p_n^{2k_n}
\cdot q_1^{2l_1+1}\cdots q_m^{2l_m+1}}$
}

\prob  A number is \dex{algebraic} iff it is the root of some polynomial with
rational coefficients. Show that the set of algebraic numbers is countable.

\prob Show that any nontrivial interval in $\rr$ has cardinality $\cc$.

\prob Show that $P(\omega)$ has cardinality $\cc$.

\prob Show that the set of all infinite subsets of $\omega$, which is
written $[\omega]^{\omega}$,
has cardinality $\cc$. \sdex{$[?omega]^{?omega}$}

\prob Show that the cardinality of $\rr \times \rr$ is $\cc$.

\prob For any set $X$ let
$[X]^{\omega}$
be the countably infinite subsets of $X$.\sdex{$[X]^{?omega}$}
Show that  $|\rr^{\omega}|=| [\rr]^{\omega}|= \cc$.
\footnote{See previous footnote.}
\com {
  $| [\rr]^{\omega}|= \cc$ requires AC. If every set of reals has
BP then there cannot be a map $g:\rr^{\omega}\mapsto \rr$ such
that $g(\langle x_n:n<\omega\rangle)=g(\langle y_n:n<\omega\rangle)$
iff $\{x_n:n<\omega\}=\{y_n:n<\omega\}$. Also can prove that
such a map gives a 1-1 map of $\omega_1$ into $\rr$, i.e.
$x_{\alpha}=f\{x_\beta:\beta<\alpha\}$ after starting with
$\{x_n:n<\omega\}$ not in the range of $f$.
}
\prob Show that the cardinality of the set of open subsets of $\rr$
is $\cc$.

\prob Show that the set of all continuous functions from $\rr$ to $\rr$
 has size $\cc$.

\prob Show that $\omega^{\omega}$  has cardinality $\cc$.

\prob Show that the set of one-to-one, onto functions from
 $\omega$ to $\omega$ has cardinality $\cc$.

\prob
Show that there is a family $\fama$ of subsets of $\qq$ such that
$|\fama | =\cc$ and
for any two distinct $s,t \in \fama$ the set  $s\cap t$ is finite.
$\fama$ is called an \dex{almost disjoint family}.

\prob Show that there is a family $\famf$ of functions from
$\omega$ to $\omega$
such that $ | \famf | = \cc$ and
for any two distinct $f,g \in \famf$ the set
$\{n \in \omega : f(n) = g(n)\}$
is finite.  These functions are called \dex{eventually different}.

\sector{Wellorderings}

 A linear order $( L, \leq)$ is a \dex{wellorder} iff for every nonempty
 $X\subseteq L$
         there exists $x \in X$ such that for every $y \in X$
        $x \leq y$   ($x$ is the minimal element of $X$).
        For an ordering $\leq$ we use $<$ to refer to the strict ordering, i.e
       $x < y$ iff $x \leq y $ and not $x \not= y$. We use $>$ to refer to the
        converse  order, i.e. $x > y$ iff $y < x$.

\bigskip

\prob Let  $( L, \leq )$ be a well ordering. Let  $(L\times L,\leq^\prime)$
 be defined  in one of    the following ways:
\par   a.  $(x,y) \leq^\prime  (u,v)$ iff $x < u$ or ($x = u$ and $y \leq v$)
\par   b.  $(x,y) \leq^\prime  (u,v)$ iff $x \leq u$ and $y \leq v$
\par   c.  $(x,y) \leq^\prime  (u,v)$ iff $max\{x,y\} < max\{u,v\} $ or
$[max\{x,y\} = max\{u,v\}$ and ($x < u$ or ($x = u$ and $y \leq v)]$.
\par\noindent   Which are well-orderings?

\prob Prove:
   Let  $( L, \leq )$ be any well-ordering and $f: L \to L$ an
    increasing function
   ($\forall x,y \in L\;(  x < y \implies f(x) < f(y))$). Then
    for every $x$ in $L$ $x \leq f(x)$.

\prob For two binary relations $R$ on $A$ and $S$ on $B$  we write
        $(A,R)\isom (B,S)$
         iff there exists a one-to-one onto map
        $f:A\to B$
        such that

\centerline{for every $x,y$ in $A$ ($xRy$ iff $f(x)Sf(y)$).}

\noindent  Such a map is called an
\dex{isomorphism}.  \sdex{$(A,R)?isom (B,S)$}
 If $(L_1,\leq_1 )$ and $(L_2,\leq_2)$ are well-orders and
$(L_1,\leq_1)\isom (L_2,\leq_2)$ then
show the isomorphism is unique. Is this true for all linear orderings?

\prob Let $(L,\leq)$ be a wellorder and for any $a \in L$ let
$L_a = \{ c \in L : c < a \}$.  Show that
 $( L, \leq )$ is not isomorphic to $( L_a, \leq )$ for any $a \in L$.

\prob (Cantor) Show that for any two wellorders  exactly one
of the following occurs:
they are isomorphic, the first is isomorphic to an initial segment
of the second, or the second is isomorphic to an initial segment
of the first.

\com{ Let $(L,\leq)$ and $(K,\leq^\prime)$ be two wellorders.
Let $$G=\{\langle a,b\rangle : ( L_a,\leq ) \isom ( K_b,\leq ^\prime )\}$$
Show that $G$ is the graph of an order preserving bijection whose domain
is all of $L$ or whose range is all of $K$.}

\prob Let $(A,\leq)$ be a linear order such that
$$\forall X\subseteq A (X\isom A \rmor (A\setminus X)\isom A ).$$
Show that $A$ is a well order or an inverse well order.

\prob Show that a linear order $(L,\leq)$ is a wellorder
 iff there does not exist an infinite sequence $x_n$ for
 $n= 0,1,2,\ldots$ with
 $x_{n+1} < x_{n}$ for every $n$.  Does this use AC?

\sector{Axiom of Choice}

\bigskip
\noindent (AC)  \dex{Axiom of Choice}:  For every family $F$ of nonempty sets
   there exists a choice function, i.e. a function $f$
   with domain F such that for every $x$ in $F$,  $f(x) \in x$.

\bigskip
\noindent (WO) \dex{Well-ordering Principle} : Every nonempty set
   can be well ordered.

\bigskip
\noindent (TL) \dex{Tuckey's Lemma}: Every family of sets
   with finite character has a maximal
   element.   A family of sets F has finite character iff for every set X,
    $X\in F$   iff for every finite $Y\subseteq X$, $Y\in F$.
\bigskip
\noindent (MP)  \dex{Maximality Principle}:  Every family of sets closed
 under the unions of chains has a maximal member.

\bigskip
\noindent (ZL) \dex{Zorn's Lemma}:  Every family of sets contains a
 maximal chain.

\par\bigskip

\prob Show that  ZL implies MP.

\prob Show that  MP implies TL.

\prob Show that  TL implies AC.

\prob (Zermelo) Show that  AC implies WO.

\com{ Given $X$  use AC to get a function
$f:P(X)\mapsto X$  such that $f(A) \in X\setminus A$
for all $A$ a proper subset of $X$.  Call a well ordering
$(L,\leq)$ respectful iff $L\subseteq X$ and for every
$a\in L\; f(L_a)=a$.
Show that the union of all respectful well orderings is a
well ordering of $X$. }

\prob  Given a nonempty family $\famf$ let $<$ be a
strict well-ordering of $\famf$.   Say that a chain $\famc\subseteq\famf$
is greedy
iff for every $a\in \famf$ if
 $$\{ b\in \famc : b < a \} \cup\{a\}$$
is a chain, then either
$a\in \famc$ or $b<a$ for every $b\in \famc$.
Show that the union of all greedy chains is a maximal chain.
Conclude that  WO implies ZL.

\prob Given a nonempty set
  $X$ let $*$ be a point
  not in  $X$ and let $Y = X \cup \{*\}$. Give $Y$ the topology
  where the open sets are $\{\emptyset,Y, X, \{ *\} \}$.
  Prove that $Y$ is a compact topological space and $X$ is
  a closed subspace of $Y$.

\prob (Kelley) The product of a family of sets is the same as
the set of all choice
functions.   Show that Tychonov's Theorem that the product of compact spaces
is compact implies the Axiom of Choice.

\sector{Ordinals}

A set $X$ is \dex{transitive} iff $\forall x \in X\,( x \subseteq X)$.
A set $\alpha$ is an \dex{ordinal} iff it is transitive and strictly
well ordered by the membership relation
(define $x\leq y$ iff $x\in y$ or $x=y$, then
$(\alpha,\leq)$ is a wellordering).  We also include the empty
set as an ordinal.
For ordinals $\alpha$ and $\beta$ we write $\alpha<\beta$ for
$\alpha \in \beta$.
The first infinite ordinal is written $\omega$. \sdex{$\omega$}
 We usually write
$$ 0 = \emptyset, 1 =\{\emptyset\}, 2 =\{\emptyset,\{\emptyset\}\}, \ldots,
n = \{ 0,1,\ldots,n-1\},\cdots, \omega =\{ 0,1,2,\ldots\} $$

\prob Show: If $\alpha$ is an ordinal then so is $\alpha + 1$.
(Remember $\alpha=\alpha\cup\{\alpha\}$.) Such
   ordinals are called \dex{successor ordinals}. Ordinals that are not
   successors are called \dex{limit ordinals}.

\prob Show: If $\alpha$ is an ordinal and $\beta < \alpha$, then $\beta$ is
 an ordinal and
$\beta \subseteq \alpha$ and $\beta=\{\gamma \in \alpha: \gamma \in \beta\}$.

\ax Axiom of Regularity.
$$ \forall x \, x \not= \emptyset \implies \exists y ( y\in x \conj
\neg \exists z\, (z \in y \conj z \in x ))$$

Another way to say this is that the binary relation
$R=\{(u,v)\in x\times x: u\in v\}$ has a minimal element, i.e.,
there exist $z$ such that for every $y\in x$ it is not the
case that $zRy$.  Note: a minimal element is not the same
as a least element.

\prob Show $ \alpha$ is an ordinal iff $\alpha$ is
   transitive and linearly ordered by the membership relation.

\prob For any ordinals $\alpha$ and $\beta$ show that
$\alpha \cap \beta  = \alpha$  or
$\alpha \cap \beta = \beta $. Show any two ordinals are comparable,
i.e., for any two distinct ordinals $\alpha$ and $\beta$ either
$\alpha\in\beta$ or $\beta\in\alpha$.

\prob The union of set $A$ of ordinals is an ordinal, and is $sup(A)$.

\prob Show that the intersection of a nonempty set $A$ of ordinals
is the least
element of $A$, written $inf(A)$.
Hence any nonempty set of  ordinals has a least element.

\prob Prove transfinite induction:  Suppose $\phi(0)$ and
$\forall \alpha\in\ord$ if $\forall \beta<\alpha\;  \phi(\beta)$,
then $\phi(\alpha)$.  Then $\forall\alpha\in\ord\; \phi(\alpha)$.

\ax Replacement Scheme Axioms.
$$\forall a\;(  [\forall x\in a \exists  ! y \,\psi (x,y)] \implies \exists b
\forall x \in a \exists y \in b \,\, \psi (x,y))$$

The formula
$\psi$ may refer to $a$ and to other sets but not to $b$.
Replacement says that for any function that is a class
the image of a set is a set. If ${\bf F}$
is a function, then for any set $a$ there exists a set $b$ such that
for every $x\in a$ there exists a $y\in b$ such that ${\bf F}(x)=y$.

\prob (von Neumann) Let $(L,\leq)$ be any well-ordering.
Show that the following is a set:
$$\{(x,\alpha): x\in L,\alpha\in\ord, \rmand (L_x,\leq)\isom \alpha\}.$$
Show that every well ordered set is isomorphic to a unique ordinal.

\bigskip
Let \ord denote the class of all ordinals. \sdex{$?ord$}
\par\noindent \dex{Transfinite Recursion}:
\par   If {\bf F} is any function defined on all sets then there exists a
            unique function {\bf G} with domain \ord such that for every
            $\alpha$ in \ord  ${\bf G}(\alpha) =
            {\bf F}( {\bf G}\res\alpha)$.

This is also referred to as a transfinite construction of ${\bf G}$.

\par\bigskip

\prob Suppose ${\bf F}:{\bf V}\to {\bf V}$, i.e., a class function.
Define g (an ordinary set function) to be a good guess iff
$\dom(g)=\alpha\in \ord$, $g(0)={\bf F}(\emptyset)$, and
$g(\beta)={\bf F}(g\res\beta)$ for every $\beta<\alpha$.
Show that if $g$ is a good guess, then $g\res\beta$ is good
guess for any $\beta<\alpha$.

\prob Show that if $g$ and $g^\prime$ are good guesses with the
same domain, then $g=g^\prime$.

\prob Show that for every $\alpha\in \ord$ there exists
a (necessarily unique) good guess $g$ with domain $\alpha$.

\prob (Fraenkel) Prove transfinite recursion.

\prob Explain the proof of WO implies ZL in terms of a
transfinite construction.

\bigskip
For an example, consider the definition of $V_{\alpha}$ for every
ordinal $\alpha$. \sdex{$V_{\alpha}$}
 Let $V_0=\emptyset$,
$V_{\alpha+1}=P(V_{\alpha})$ (power set) for successor
ordinals, and for
limit ordinals $V_{\lambda}=\cup_{\beta<\lambda} V_{\beta}$.
Thus if we define ${\bf F}$ as follows:
$${\bf F}(x)=\left\{
\begin{array}{ll}
    P(x(\alpha))  & \mbox{if $x$ is a function
                     with domain $\alpha+1\in\ord$} \\
    \bigcup_{\alpha<\lambda} x(\alpha)
                  & \mbox{if $x$ is a function
                     with domain limit ordinal $\lambda$} \\
    P(\emptyset)  & \mbox{if $x=\emptyset$} \\
    -\pi          & \mbox{otherwise}\\
\end{array}
\right.$$
then ${\bf G}(\alpha)=V_\alpha$.

\prob Show that if $\alpha\leq\beta$, then $V_\alpha\subseteq V_{\beta}$
and if $\alpha<\beta$, then $V_{\alpha}\in V_\beta$.
Show that each $V_{\alpha}$ is transitive.

\prob Show that every set is included in a transitive set.
Show that for every transitive set $x$ there exists an ordinal
$\alpha$ such that $x\in V_{\alpha}$, i.e.,
${\bf V}=\cup_{\alpha\in \ord}V_{\alpha}$.

\par\bigskip
  \dex{Ordinal arithmetic}: ($\alpha,\beta$ are ordinals, and
$\lambda$ is a limit ordinal.)

\medskip
Addition:

$\alpha + 0 = \alpha$

$\alpha +(\beta+1)=(\alpha+\beta)+1$

$\alpha+\lambda=sup \{ \alpha + \beta : \beta < \lambda \}$

\medskip
Multiplication:

$\alpha   0 =  0$

$ \alpha  ( \beta + 1)  =  (\alpha   \beta) + \alpha$

$\alpha\lambda=  sup \{ \alpha   \beta : \beta < \lambda \}$

\medskip
Exponentiation:

$ \alpha^0=1$

$\alpha^{\beta+1}=\alpha^\beta\alpha$

$\alpha^{\lambda}=sup \{ \alpha^{\beta}: \beta < \lambda \}$

\medskip
 So for example the addition function
$+ : \ord \times \ord \mapsto
\ord$ exists by transfinite recursion.
For each $\alpha \in$ \ord define a function ${\bf F}_{\alpha}$
on all sets as follows:

${\bf F}_{\alpha}(g)=g(\beta)+1$ if
 $g$ is a map with domain an ordinal $\beta+1$,

${\bf F}_{\alpha}(g)= sup\{g(\gamma): \gamma < \lambda\}$ if
 $g$ is a map with domain a limit  ordinal $\lambda$, and

${\bf F}_{\alpha}(g)=\alpha$ otherwise.

\noindent Hence for each $\alpha$ we have a unique
${\bf G}_{\alpha}: \ord \mapsto \ord$
which will exactly be
${\bf G}_{\alpha}(\beta) = \alpha + \beta$. Since each ${\bf G}_{\alpha}$ is
 unique we have defined  the function $+$ on all pairs of ordinals.

Note that more intuitively $\alpha+\beta$ is the unique ordinal isomorphic
to the well-order $(\{0\}\times\alpha)\cup (\{1\}\times\beta)$ ordered
lexicographically.  Similarly $\alpha\cdot\beta$ is the unique ordinal
isomorphic to the well-order $\beta\times\alpha$ ordered
lexicographically.  Exponentiation is much harder to describe.

\bigskip

\prob Show that $\alpha + (\beta +\gamma)=(\alpha + \beta) +\gamma$.

\prob Assume $\alpha$, $\beta$, and $\gamma$ are ordinals,
 which of the following are always true?
\begin{eqnarray*}
    \alpha + \beta & = & \beta + \alpha  \\
    \alpha + ( \beta +\gamma ) & = & ( \alpha + \beta ) +\gamma  \\
    \alpha   \beta & = & \beta   \alpha                            \\
    \alpha   ( \beta   \gamma )  & = & ( \alpha   \beta )
             \gamma  \\
    \alpha   ( \beta +\gamma ) & = & ( \alpha   \beta ) + ( \alpha
        \gamma )   \\
    ( \alpha + \beta )  \gamma & = & ( \alpha  \gamma ) +
       ( \beta  \gamma )     \\
     \alpha^{ \beta   \gamma} & = & \bigl(\alpha^{\beta}\bigr)^{\gamma}  \\
     \alpha^{\beta}  \alpha^{\gamma} & = &  \alpha^ {\beta +\gamma}\\
    (\alpha\beta)^\gamma=\alpha^\gamma\beta^\gamma\\
     \alpha + \beta = \alpha + \gamma & \implies & \beta =\gamma \\
     \beta < \gamma & \iff & \alpha + \beta < \alpha + \gamma \\
     \beta + \alpha = \gamma + \alpha & \implies & \beta =\gamma \\
  (\alpha>0 \;\conj\;  \alpha   \beta = \alpha   \gamma)
                      & \implies & \beta =\gamma \\
  (\alpha>1\;\conj\;   {\beta <\gamma})
                      & \implies & \alpha^{\beta} < \alpha^{\gamma} \\
\end{eqnarray*}

\prob For any ordinals $\alpha$ and $\beta>0$ show there exists
unique ordinals
 $\gamma$ and $\delta$ such that
     $\alpha =  \beta   \gamma   +  \delta$ and $\delta<\beta$.

\prob Is the previous problem true for  $\alpha=\gamma\beta + \delta$?

\prob Show for any $\beta > 0$ there exists unique ordinals
    $\gamma_1,\ldots,\gamma_n, d_1,\ldots,d_n$  such that
    $\gamma_1 >\gamma_2 > \cdots >\gamma_n$ ;
    $ 0 <  d_1,d_2,\ldots,d_n < \omega$ and
$$ \beta = \omega^{\gamma_1}d_1 + \omega^{ \gamma_2}d_2 + \ldots +
       \omega^{ \gamma_n}d_n $$
This is called \dex{Cantor normal form}.

\prob Sort the following set of five ordinals:

$$\begin{array}{ccc} (\omega^{\omega})(\omega+\omega) &
(\omega+\omega)   ( \omega^{\omega} ) &
 \omega^{\omega}   \omega + \omega^{\omega}   \omega\\
\omega  \omega^{\omega} +\omega   \omega^{\omega} &
\omega^{\omega}  \omega+\omega  \omega^{\omega}
\end{array} $$

\prob An ordinal  $\alpha$  is \dex{indecomposable} iff it satisfies any
 of the following.

 a) $\exists \beta \,\,\alpha = \omega^{\beta} $

 b) $\forall \beta \forall \gamma$ if $\alpha= \beta +\gamma$ then
        $\alpha=\beta$ or $\alpha =\gamma$

 c) $\forall X \subseteq \alpha \,\,\,[( X,<) \isom (\alpha,<)$ or
         $(\alpha\setminus X,<) \isom (\alpha,<)]$

 d) $\forall \beta<\alpha\;\;\beta+\alpha=\alpha$

Show they are all equivalent.

\prob (Goodstein) The complete expansion of a positive integer in a base
 is gotten by writing
everything possible including exponents and exponents of exponents etc. as
powers of the base.  For example the number 36 written in complete base
two is: $$ 2^{( {2^2} + 1)}+2^2$$
The same number in complete base 3 is:
$$ 3^3 + 3^2$$
Let $a_n$ be a sequence described as follows.  Given $a_k$ calculate
$a_{k+1}$ by
writing $a_k$ in base k then substitute k+1 for every k, then subtract one.
For example:
$$
a_2 =36=2^{( {2^2} + 1)}+2^2  \implies  a_3 =
3^{\Bigl( {3^3} + 1 \Bigr)}+3^3 - 1
 = 2.2876\ldots \times 10^{13}
$$
or for example if $a_6 = 6^4 +  2 \cdot  6^3 = 1728$, then
$$
 a_7 =  ( 7^4   + 2 \cdot  7^3) -1
= 7^4 +7^3 +6\cdot  7^2 + 6\cdot 7^1 + 6 = 3086
$$
Show that given any pair of positive integers n and m, if we let
$a_n = m$ then for some $k > n$ we get $a_k = 0$.

\prob Let $S$ be a countable set of ordinals.  Show that
$$\{\beta:\exists \langle \alpha_n:n\in\omega\rangle\in S^{\omega}
\;\;\beta=\Sigma_{n<\omega}\alpha_n \}$$
is countable.
($\Sigma_{n<\omega}\alpha_n=sup\{\Sigma_{n<m}\alpha_n:m<\omega\}$)

\sector{Cardinal Arithmetic}

An ordinal $\kappa$ is a \dex{cardinal} iff for every $\alpha < \kappa$,
$|\alpha| < | \kappa|$.
The cardinality of a set A is the least cardinal $\kappa$,
$|A|=|\kappa|$ and we write $|A|=\kappa$.
The $\alpha^{th}$ uncountable cardinal is written either
$\aleph_{\alpha}$
or $\omega_{\alpha}$.  \sdex{$\aleph_{\alpha}$}
\sdex{$\omega_{\alpha}$}

\bigskip

\prob Show that for any ordinal $\alpha$ the cardinal
$\aleph_\alpha$ exists.

\prob Is there a cardinal such that $\aleph_{\alpha}=\alpha$ ?

\bigskip

For cardinals $\kappa$ and $\gamma$  we define  $\kappa \gamma$
   to be the cardinality of the
cross product and $\kappa +\gamma$ to be the cardinality of the union
 of $A$ and $B$ where $A$ and $B$ are disjoint and
 $|A|=\kappa$ and $|B|=\gamma$.  \sdex{$\kappa \gamma$}
 \sdex{$\kappa +\gamma$}

\bigskip

\prob Let $\kappa$ be an infinite cardinal.
Define the \dex{lexicographical order} $\leq_l$ on  $\kappa \times \kappa$
by $(x,y) \leq_l  (u,v)$ iff $x<u$ or $(x=u$ and $y\leq v)$.
Define $ \leq^{\prime}$ on $\kappa \times \kappa$ by
$(x,y) \leq^{\prime}  (u,v)$ iff
$max\{x,y\} < max\{u,v\}$  or
$[(max\{x,y\} = max\{u,v\}$ and $(x,y)\leq_l (u,v)]$.
Show that  $( \kappa, \leq ) \isom ( \kappa \times \kappa , \leq^{\prime})$.

\prob Show that for infinite cardinals $\kappa$ and $\gamma$,
  $\kappa +\gamma = \kappa \gamma  = max \{ \kappa,\gamma\}$.

\prob Show that for any infinite cardinal $\kappa$ the union of $\kappa$
 many sets of cardinality $\kappa$ has cardinality $\kappa$.

\bigskip

The \dex{cofinality} of an infinite limit ordinal $\beta$,
$cf(\beta)$,  \sdex{$cf(?beta)$}
 is the least $\alpha\leq\beta$
 such that there
is a map $f:\alpha\to\beta$ whose range is unbounded in $\beta$.

\bigskip

\prob For $\lambda$ a limit ordinal show the following
are all equivalent:

a. $\alpha$ is the minimum ordinal such that
$\exists X \subseteq \lambda$ unbounded in $\lambda$
such that $(X,\leq) \isom (\alpha,\leq)$.

b. $\alpha$ is the minimum ordinal such that
$\exists f:\alpha\to\lambda$ such that
f is one-to-one,  order preserving, and the range of $f$ is
cofinal (unbounded) in $\lambda$.

c.  $cf(\lambda)=\alpha$

\prob For $\kappa$ an infinite cardinal show that
$cf(\kappa)=\alpha$ iff $\alpha$ is the minimum cardinal such that
$\kappa$ is the union of $\alpha$ many sets of cardinality less
than $\kappa$.

\prob Let $\alpha$ and $\beta$ be limit ordinals and suppose $f:\alpha\to
\beta$ is strictly increasing and cofinal in $\beta$.
Show $cf(\alpha)=cf(\beta)$.

\prob $\kappa$ is \dex{regular} iff $cf(\kappa)=\kappa$.
$\kappa$ is \dex{singular} iff $cf(\kappa)<\kappa$.
Show that for any limit ordinal $\beta$,
$cf(\beta)$ is a regular cardinal.

\prob
$\kappa^+$ is the least cardinal greater than $\kappa$. \sdex{$\kappa^+$}
Show that for any infinite cardinal $\kappa$ ,  $\kappa^+$ is a
regular cardinal.

\prob For $\alpha$ a limit ordinal show that
$cf(\aleph_\alpha)=cf(\alpha)$.

\bigskip

$\kappa^{\gamma}$  is the set of all functions from $\gamma$
   to $\kappa$, but we often use it to denote its own cardinality.
   \sdex{$\kappa^{\gamma}$}
  $\kappa^{<\gamma}$ is  $\cup\{\kappa^{\alpha}: \alpha < \gamma\}$,
  but we often use it to denote its own cardinality, i.e.,
  $\kappa^{<\gamma}=|\kappa^{<\gamma}|$. \sdex{$\kappa^{<\gamma}$}
{\bf Note in the section all exponentiation is cardinal exponentiation.}

\bigskip

\prob Show that $|\kappa^{<\omega}|=|[\kappa]^{<\omega}|=\kappa$ for any
infinite cardinal $ \kappa$. $[\kappa]^{<\omega}$ is the set
of finite subsets of $\kappa$. \sdex{$[\kappa]^{<\omega}$}

\prob Show that $\aleph_{\omega} <  \aleph_{\omega}^{\omega}$. Show
that for any cardinal $\kappa$, $\kappa < \kappa^{cf(\kappa)}$.

\prob (K\"{o}nig) Show that $cf(2^{\kappa})>\kappa$.

\prob Show $(\forall n \in \omega \,\, 2^{\aleph_n}=
\aleph_{n+1}) \implies 2^{\aleph_{\omega}}=
(\aleph_{\omega})^{\omega}$.

\prob Show that $({2^{<\kappa}})^{cf(\kappa)}=2^{\kappa}$.

\prob Show $(\forall n \in \omega \,\, 2^{\aleph_n}=
\aleph_{\omega+17})\implies 2^{\aleph_{\omega}}=
\aleph_{\omega+17}$.

\prob Let $\kappa$ be the least cardinal such that
$2^{\kappa}> 2^{\omega}$. Show that $\kappa$ is regular.

\prob Prove that for every infinite regular cardinal $\kappa$, there is
a cardinal $\lambda$ such that $\aleph_{\lambda}=\lambda$ and
$\lambda$ has cofinality $\kappa$.

\prob Show that $cf(2^{<\kappa}) = cf(\kappa)$ or
 $cf(2^{<\kappa}) >\kappa$.

\prob Show that if $\omega \leq \lambda \leq \kappa$ then
$(\kappa^{+})^{\lambda}= max\{\kappa^{\lambda},\kappa^{+}\}$.

\prob For any set $X$ \dex{Hartog's ordinal} $h(X)$ is defined by:
$$ h(X) = sup\{ \alpha \in ORD: \exists {\rm \,onto\,} f: X \to \alpha\}$$
Show without AC that $h$ is well defined and with AC $h(X)=|X|^+$.
Show that AC is equivalent to the statement that for every two sets $X$ and
$Y$ either $|X|\leq|Y|$ or $|Y|\leq|X|$.

\prob Given sets $A_{\alpha} \subseteq \kappa$ for each $\alpha <\kappa$ each
of cardinality $\kappa$, show there exists $X\subseteq\kappa$ such that
$$\forall\alpha <\kappa\,\,\,\,
|A_{\alpha}\cap X|=|A_{\alpha}\setminus X|=\kappa$$

\prob Show there exists $X\subseteq \Bbb R$ which contains the rationals
and the only
automorphism of $(X,\leq)$ is the identity, i.e., any order preserving
bijection from $X$ to $X$ is the identity.

\prob Show there exists $X\subseteq \Bbb R^2$ such that for every $x\in X$
and positive $ r\in \Bbb R$ there is a unique $y\in X$ with
$d(x,y)=r$ where $d$ is Euclidean distance.

\prob (Sierpi\'{n}ski) Show that there exists
$X\subseteq \Bbb R^2$ such that for every line
$L$ in the plane, $|L\cap X|=2$.

\prob (Jech) Without using AC show that $\omega_2$ is not the countable
union of countable sets.

\sector{First Order Logic and the Compactness Theorem}

\begin{center}
  Syntax
\end{center}

We begin with the syntax of first order logic. The \dex{logical symbols}
are $\disj,\neg,\exists,=$ and for each $n\in\omega$ a variable symbol
$x_n$.  There are also grammatical symbols such as parentheses and commas
that we use to parse things correctly but have no meaning.
For clarity we usually use $x,y,z,u,v,$ etc. to refer to
 arbitrary variables.
The \dex{nonlogical symbols} consist of a given set $L$ that
may include operation
symbols, predicate symbols, and constant symbols.  The case where $L$ is
empty is referred to as the \dex{language of pure equality}.
Each symbol $s\in L$ has a nonnegative integer $\#(s)$ called its
\dex{arity}
assigned to it.
If $\#(s)=0$, then $s$ is a constant symbol. If $f$ is an operation symbol
and $\#(f)=n$ then $f$ is an $n$-ary operation symbol. Similarly if $R$ is
a predicate symbol and $\#(R)=n$ then $R$ is an $n$-ary predicate symbol.
In addition we always have that ``$=$'' is a logical binary predicate symbol.

For the \dex{theory of groups} the appropriate
language is $L=\{e,\cdot,^{-1}\}$ where ``$e$'' is a constant symbol,
so $\#(e)=0$,
``$\cdot$'' is a binary operation symbol, so  $\#(\cdot)=2$, and
``$^{-1}$''
 is a unary operation symbol, so $\#( ^{-1})=1$.
   For the \dex{theory of partially ordered sets} we have that
$L=\{\leq\}$ where $\leq$ is a binary relation symbol, so $\#(\leq)=2$.

Our next goal is to define what it means to be a formula
 of first order logic.
Let $L$ be a fixed language.   An expression is a finite string of
 symbols that are either logical symbols or symbols from $L$.

The set of \dex{terms} of $L$
is the smallest set of expressions that contain the variables and
constant symbols of $L$ (if any), and is closed under the formation rule:
if $t_1,t_2,\ldots,t_n$  are terms of $L$ and  $f$ is an $n$-ary
 operation symbol of $L$, then $t=f(t_1,t_2,\ldots,t_n)$ is a term of $L$.
If $L$ has no function symbols then the only terms of $L$ are the
variables and constant symbols.  So for example if $c$ is a constant symbol,
$f$ is a 3-ary operation symbol, $g$ is a binary operation symbol, and
$h$ is a unary operation symbol, then
$$h(f(g(x,h(y)),y,h(c)))$$
is a term.

The set of \dex{atomic formulas} of $L$ is the set of all expressions of the
form $R( t_1,t_2,\ldots,t_n)$, where  $t_1,t_2,\ldots,t_n$ are terms
of $L$ and $R$ is a $n$-ary predicate symbol of $L$.  Since we
always have equality as
a binary relation we always have atomic formulas of the form
 $t_1=t_2$.

The set of \dex{formulas} of $L$ is
the smallest set of expressions that includes
the atomic formulas and is closed under the formation rule:
if $\theta$ and $\psi$ are $L$ formulas and $x$ is any variable, then
\begin{itemize}
\item $(\theta \disj \psi)$,
\item $\neg\theta$, and
\item$\exists x \,\theta$
\end{itemize}
are $L$
formulas.  We think of other logical connectives as being abbreviations,
e.g.,
\begin{itemize}
\item $(\theta \conj\psi)$ abbreviates $\neg(\neg\theta \disj \neg\psi)$,
\item $(\theta \implies\psi)$ abbreviates $(\neg\theta \disj \psi)$,
\item $\forall x \,\theta $ abbreviates $\neg\exists x \,\neg\theta$,
and so forth.
\end{itemize}

We often add and sometimes drop parentheses to improve readability.
Also we write $x\not=y$ for the formally correct but harder to
read $\neg x=y$.

It is common practice to write symbols not only in \dex{prefix} form as above
but also in \dex{postfix} and \dex{infix} forms.
For example in our example of
group theory instead of writing the term $\cdot(x,y)$ we usually write it
in infix form $x \cdot y$, and $^{-1}(x)$ is usually written
in postfix form $x^{-1}$.  Similarly in the language of partially ordered sets
 we usually write
 $x\leq y$ instead of the prefix form $\leq(x,y)$.  Binary relations
  such as partial
orders and equivalence relations are most often written in infix form.
We regard the more natural forms we write as abbreviations of the more
formally correct prefix notation.

Next we want to describe the syntactical concept of \dex{substitution}.
To do so we must first describe what it means for an occurrence of a
variable $x$ in a formula $\theta$ to be \dex{free}. If an occurrence of a
variable $x$ in a formula $\theta$ is not free it is said to be \dex{bound}.
Example:
$$ (\exists x\,\, x = y \disj x = f(y) )$$
Both occurrences of $y$ are free, the first two occurrences of $x$ are bound,
and the last occurrence of $x$ is free. In the formula:
$$ \exists x\,\, (x = y \disj x = f(y) )$$
all three occurrences of $x$ are bound.

 Formally we proceed as follows.  All occurrences of variables
in an atomic formula are free.  The free occurrences of $x$ in $\neg\theta$
are exactly the free occurrences of $x$ in $\theta$.
 The free occurrences
of $x$ in $(\theta\disj\psi)$ are exactly the free occurrences of $x$ in
$\theta$ and in $\psi$. If $x$ and $y$ are distinct variables, then the free
occurrences of $x$ in $\exists y\,\,\theta$ are exactly the free occurrences
of $x$ in $\theta$.  And finally no occurrence of $x$ in $\exists x\,\,\theta$
is free. This gives the inductive definition of free and bound variables.

We show that $x$ might occur freely in $\theta$ by
writing $\theta(x)$. If $c$ is a constant symbol the formula $\theta(c)$
is gotten by substituting $c$ for all free occurrences (if any) of
$x$ in $\theta$.  For example: if $\theta(x)$ is
$$\exists y\; (y=x\conj\forall x\; x=y),$$
then $\theta(c)$ is
$$\exists y\; (y=c\conj\forall x\;x=y).$$

We usually write $\theta( x_1,x_2,\ldots,x_n)$ to indicate that the free
variables of $\theta$ are amongst the $x_1,x_2,\ldots,x_n$.
A formula is called a sentence if no variable occurs freely in it.

\begin{center}
  Semantics
\end{center}

 Our next goal is to describe the semantics of first order logic.
 A \dex{structure} $\moda$ for the language $L$ is a pair consisting
 of a set $A$
 called the universe of $\moda$ and an assignment or interpretation
 function from the nonlogical symbols
 of $L$ to individuals, relations, and functions on $A$.  Thus
\begin{itemize}
\item for each
 constant symbol $c$ in $L$ we have an assignment $c^{\moda}\in  A$,

\item for each $n$-ary operation symbol $f$ in $L$ we have a function
 $f^{\moda}: A^n\to A$, and

\item for each $n$-ary predicate symbol $R$ we
 have a relation $R^{\moda}\subseteq A^n$.
\end{itemize}
The symbol $=$ is always
interpreted
 as the binary relation of equality,  which is why we consider it a
 logical symbol, i.e., for any structure $\moda$ we have
 $=^{\moda}$ is $\{(x,x):x\in A\}$.   We use the word structure and
 \dex{model} interchangeably.

For example, suppose $L$ is the language of group theory.
One structure for
this theory is
$$\modq=({\Bbb Q},+,-x,0)$$
where
\begin{itemize}
\item the universe is the rationals,
\item $\cdot^{\modq}$ is ordinary addition of rationals,
\item $^{-1^{\modq}}$ is the function which takes each rational
      $r$ to  $-r$, and
\item $e^{\modq}=0$.
\end{itemize}
Another structure in this language is
$$\modr=({\Bbb R}^+,\times,{1 \over x},1)$$
where
\begin{itemize}
\item the universe is the set of positive real numbers,
\item $\cdot^{\modr}$ is multiplication $\times$,
\item $^{-1^{\modr}}$ is the function which takes $x$ to $1 \over x$, and
\item $e^{\modr}=1$.
\end{itemize}
Another example is the  group $S_n$ of permutations.
Here $\cdot^{S_n}$ is composition of functions,
$^{-1^{S_n}}$ is the functional
which takes each permutation to its inverse, and $e^{S_n}$ is the identity
permutation.   Of course there are many examples of structures in this
language which are not groups.

For another example, the language of partially ordered sets is
$L=\{\leq\}$ where $\leq$ is a binary relation symbol.  The
following are all $L$-structures which happen to be partial orders:
\begin{itemize}
\item $({\Bbb R},\{(x,y)\in {\Bbb R}^2: x\leq y\})$,
\item $({\Bbb Q},\{(x,y)\in {\Bbb Q}^2: x\geq y\})$, and
\item $({\Bbb N},\{(x,y)\in {\Bbb N}^2: x\mbox{ divides }y\})$.
\end{itemize}
For any nonempty set
$A$ and  $R\subseteq A^2$, $(A,R)$ is an $L$-structure.  If in addition
the relation $R$ is transitive, reflexive, and antisymmetric, then
$(A,R)$ is a partial order.

\begin{center}
  ${\moda\models\theta}$
\end{center}

Next we define what it means for an $L$ structure $\moda$ to model or
satisfy
an $L$ sentence $\theta$ (written ${\moda\models\theta}$).
\sdex{$?moda?models?theta$}
For example,
$$({\Bbb Q},+,0)\models \forall x \exists y \; x\cdot y=e,$$
because for all $p\in {\Bbb Q}$ there exists $q\in {\Bbb Q}$ such that
$p+q=0$.

Usually it is not the case that every element of a model has a
constant symbol which names it.  But suppose this just happened to
be the case.  Let's suppose that for ever $a\in A$ there is a
constant symbol $c_a$ in the language $L$ so that
$c_a^\moda=c$.
The interpretation function can be extended to the variable free terms
of $L$ by the rule:
$$ (f( t_1,t_2,\ldots,t_n))^{\moda} =
f^{\moda}( t_1^{\moda},t_2^{\moda},\ldots,t_n^{\moda})$$
Hence for each variable free term $t$ we get an interpretation
$t^{\moda}\in A$.  For example, if $L=\{S,c\}$ where $S$ is
a unary operation symbol and $c$ is a constant symbol, and
${\goth Z}$ is the $L$-structure with universe ${\Bbb Z}$ and
where $S^{\goth Z}(x)=x+1$ and $c^{\goth Z}=0$, then
$S(S(S(S(c))))^{\goth Z}=4$.

Our definition of $\models$ is by induction on the \dex{logical complexity}
 of the sentence $\theta$, i.e. the number of logical symbols in $\theta$.

\begin{enumerate}
\item $\moda\models R(t_1,\ldots,t_n)$ iff
 $( t_1^{\moda},\ldots,t_n^{\moda})\in R^{\moda}$.

\item $\moda\models\neg\theta$ iff not $\moda\models\theta$.

\item $\moda\models(\theta\disj\psi)$ iff $\moda\models\theta$
or $\moda\models\psi$.

\item $\moda\models\exists x\,\theta(x)$ iff there exists a $b$ in
the universe $A$ such that  $\moda\models\theta(c_b)$.
\end{enumerate}

Now we would like to define $\moda\models\theta$ for arbitrary
languages $L$ and $L$-structures $\moda$.
Let ${L_{\moda}=L\cup\{c_a: a\in A\}}$
\sdex{$L_{?moda}=L?cup?{c_a: a?in A?}$}
 where each ${c_a}$ is a new constant symbol.
 Let ${(\moda,a)_{a\in A}}$ be the $L_{\moda}$
\sdex{$(?moda,a)_{a?in A}$}
structure gotten by augmenting the structure $\moda$ by interpreting
each symbol $c_a$ as the element $a$.

If $\theta$ is an $L$-sentence and
$\moda$ is and $L$-structure, then we define {$\moda\models\theta$}
iff $\theta$ is true in the augmented structure, i.e., ${(\moda,a)_{a\in A}}$.

If $L_1 \subseteq L_2$ and $\moda$ is a $L_2$ structure, then the reduct of
$\moda$
to $L_1$, written $\moda\res{L_1}$, is the $L_1$ structure with the
same universe   \sdex{$?moda?res{L_1}$}
as $\moda$ and same relations, operations, and constants as $\moda$
for the symbols of $L_1$.

\begin{lemma}
Let $L_1\subseteq L_2$ and $\moda$ be an $L_2$ structure.  Then
for any $\theta$ an $L_1$ sentence,
$$\moda\models\theta \mbox{ iff }
\moda\res{L_1}\models\theta$$
\end{lemma}
\proof
We prove by induction on the number of logical symbols in the sentence
that for any $L_{1\moda}$ sentence $\theta$:
$$(\moda,a)_{a\in A}\models\theta \mbox{ iff }
(\moda,a)_{a\in A}\res{L_{1\moda}} \models\theta$$
Let $\moda_2=(\moda,a)_{a\in A}$ and
    $\moda_1=(\moda,a)_{a\in A}\res{L_{1\moda}}$.

\medskip
Atomic sentences: By induction  on the
size of the term, for any $L_{1\moda}$ variable free term $t$
we have that $t^{\moda_1}=t^{\moda_2}$.

For any $n$-ary relation
symbol $R$ in $L_1$ we have $R^{\moda_1}=R^{\moda_2}$ (since
$\moda_1$ is a reduct of $\moda_2$).  Hence for any atomic
$L_{1\moda}$-sentence $R(t_1,\ldots,t_n)$,
$$\moda_1\models R(t_1,\ldots,t_n)$$
\centerline{iff}
$$\langle t_1^{\moda_1},\ldots t_n^{\moda_1}\rangle \in R^{\moda_1}$$
\centerline{iff}
$$\langle t_1^{\moda_2},\ldots t_n^{\moda_2}\rangle \in R^{\moda_2}$$
\centerline{iff}
$$\moda_2\models R(t_1,\ldots,t_n).$$

\medskip
Negation:
$$\moda_1\models \neg\theta$$
\centerline{iff}
$$\mbox{ not } \moda_1\models \theta$$
\centerline{iff (by induction)}
$$\mbox{ not }\moda_2\models \theta$$
\centerline{iff}
$$\moda_2\models \neg\theta.$$

\medskip
Disjunction:
$$\moda_1\models (\theta\disj\rho)$$
\centerline{iff}
$$\moda_1\models \theta \rmor \moda_1\models\rho$$
\centerline{iff (by induction)}
$$\moda_2\models \theta \rmor \moda_2\models\rho$$
\centerline{iff}
$$\moda_2\models (\theta\disj\rho).$$

\medskip
Existential quantifier:
$$\moda_1\models \exists x\theta(x)$$
\centerline{ iff}
$$\mbox{ there exists $a\in A$ such that }\moda_1\models \theta(c_a)$$
\centerline{ iff}
$$\mbox{ there exists $a\in A$ such that }\moda_2\models \theta(c_a)$$
\centerline{iff}
$$\moda_2\models \exists x\theta(x).$$

\qed

\begin{center}
   Compactness Theorem
\end{center}

The compactness theorem (for countable languages)
was proved by Kurt G\"{o}del in 1930.  Malcev extended it to uncountable
languages in 1936.  The proof we give here was found by Henkin in 1949.

We say that a set of $L$ sentences $\Sigma$ is \dex{finitely satisfiable}
iff every finite subset of $\Sigma$ has a model.  $\Sigma$ is
\dex{complete} iff for every $L$ sentence $\theta$ either
$\theta$ is in $\Sigma$ or $\neg\theta$ is in $\Sigma$.

\begin{lemma}
 For every finitely satisfiable set of $L$ sentences $\Sigma$
 there is a complete finitely satisfiable set of $L$ sentences
 $\Sigma^\prime\supseteq \Sigma$.
\end{lemma}
\proof
Let $B=\{Q:Q\supseteq \Sigma$ is finitely satisfiable$\}$.
$B$ is closed under unions of chains, because if $C\subseteq B$ is
a chain, and $F\subseteq \cup C$ is finite then there exists
$Q\in C$ with $F\subseteq Q$, hence $F$ has a model.  By the
maximality principal, there exists $\Sigma^\prime\in B$ maximal.
But for every $L$ sentence $\theta$ either $\Sigma^\prime\cup\{\theta\}$
is finitely satisfiable or $\Sigma^\prime\cup\{\neg\theta\}$ is
finitely satisfiable.  Otherwise there exists finite
$F_0,F_1\subseteq\Sigma^\prime$ such that
$F_0\cup\{\theta\}$  has no model and $F_1\cup\{\neg\theta\}$
has no model.  But  $F_0\cup F_1$ has a model $\moda$ since
$\Sigma$ is finitely satisfiable.
Either $\moda\models\theta$ or $\moda\models\neg\theta$.
This is a contradiction.
\qed

\begin{lemma} If $\Sigma$ is a finitely satisfiable set of $L$ sentences,
and $\theta(x)$ is an $L$ formula with
one free variable $x$,
and $c$ a new constant symbol (not in $L$), then
$\Sigma \cup \{(\exists x \,\theta (x)) \implies \theta(c)\}$
is finitely satisfiable in the language $L\cup\{c\}$.
\end{lemma}
\proof
This new sentence is called a \dex{Henkin sentence} and $c$ is called the
\dex{Henkin constant}.  Suppose it is not finitely satisfiable.  Then
there exists $F\subseteq\Sigma$ finite such that
$F\cup\{(\exists x \,\theta (x)) \implies \theta(c)\}$ has no
model.  Let $\moda$ be an $L$-structure modeling $F$.
Since the constant $c$ is not in the language $L$ we are free to
interpret it any way we like.  If $\moda\models\exists x \,\theta (x)$
choose $c\in A$ so that $(\moda,c)\models\theta(c)$, otherwise
choose $c\in A$ arbitrarily.  In either case
$(\moda,c)$ models $F$ and the Henkin sentence.
\qed

\par\medskip
We say that a set of $L$ sentences $\Sigma$ is \dex{Henkin}
iff for every $L$ formula  $\theta(x)$ with one free variable $x$,
there is a constant symbol   $c$ in $L$ such
that $(\exists x \,\theta (x)) \implies \theta(c) \in \Sigma$.

\begin{lemma} If $\Sigma$ is a finitely satisfiable set of $L$ sentences,
 then there exists
 $\Sigma^\prime\supseteq \Sigma$
with $L^\prime\supseteq L$ and $\Sigma^\prime$ a finitely satisfiable
Henkin set of $L^\prime$ sentences.
\end{lemma}
\proof
For any set of $\Sigma$ of $L$ sentences, let
$$\Sigma^*=\Sigma\cup\{
(\exists x \,\theta (x)) \implies \theta(c_{\theta}):
 \theta(x)\mbox{ an }L\mbox{ formula with one free variable}\}$$
The language of $\Sigma^*$ contains a new constant symbol $c_{\theta}$
for each $L$ formula $\theta(x)$.  $\Sigma^*$ is finitely satisfiable,
since any finite subset of it is contained
in a set of the form
 $$F\cup\{(\exists x \,\theta_1 (x)) \implies \theta(c_{\theta_1}),\ldots,
            (\exists x \,\theta_n (x)) \implies \theta(c_{\theta_n})\}$$
where  $F\subseteq\Sigma$ is finite.  To prove this set has a model
use induction on $n$ and note that
from the point of view of
$$\Sigma\cup\{(\exists x \,\theta_1 (x))
\implies \theta(c_{\theta_1}),\ldots,
(\exists x \,\theta_{n-1} (x)) \implies \theta(c_{\theta_{n-1}})\}$$
$c_{\theta_n}$ is a new constant symbol, so we can apply the last lemma.

Now let $\Sigma_0=\Sigma$ and let
$\Sigma_{m+1}=\Sigma_m^*$.  Then
$$\Sigma^\prime=\cup_{m<\omega}\Sigma_m$$
is Henkin.  It is also finitely satisfiable, since it is the union
of a chain of finitely satisfiable sets.
\qed

\par\medskip
If $\Sigma$ is a set of $L$ sentences, then the \dex{canonical structure}
$\moda$ built
from $\Sigma$ is the following.
\begin{itemize}
\item Let $X$ be the set of all variable free terms of $L$.
\item For $t_1,t_2\in X$ define $t_1\sim t_2$ iff $(t_1 = t_2) \in \Sigma$.
\item Assuming that $\sim$ is an equivalence relation let
$[t]$ be the equivalence
class of $t\in X$.
\item The universe of the canonical model $\moda$ is the set of
equivalence classes of $\sim$.
 \item For any constant symbol $c$ we define
$$c^{\moda}=[c].$$
\item For any $n$-ary operation symbol $f$ we define
$$f^{\moda}([t_1],[t_2],\ldots,[t_n])=[f(t_1,t_2,\ldots,t_n)].$$
\item For any $n$-ary relation symbol $R$ we define
$$([t_1],[t_2],\ldots,[t_n])\in R^{\moda} \rmiff
 R(t_1,t_2,\ldots,t_n)\in \Sigma.$$
\end{itemize}

\begin{lemma} If $\Sigma$ is a finitely satisfiable complete Henkin set of
$L$ sentences, then the canonical
model $\moda$ built from $\Sigma$ is well defined and for every $L$
sentence $\theta$,
$$\moda\models \theta \rmiff \theta \in \Sigma.$$
\end{lemma}
\proof
First we show that $\sim$ is an equivalence relation.
Suppose $t,t_1,t_2,t_3$ are variable free terms.

\smallskip
$t\sim t$:
If $t=t\notin \Sigma$ then,
since $\Sigma$ is complete we have that $\neg t=t\in \Sigma$.
But clearly $\neg t=t$ has no models and so $\Sigma$ is not finitely
satisfiable.

\smallskip
$t_1\sim t_2$ implies $t_2\sim t_1$:
If not, by completeness of $\Sigma$ we must have
that $t_1=t_2$ and $\neg t_2=t_1$ are both in $\Sigma$. But then
$\Sigma$ is not finitely satisfiable.

\smallskip
($t_1\sim t_2$ and $t_2\sim t_3$) implies $t_1\sim t_3$:
If not, by completeness of $\Sigma$ we must have
that $t_1=t_2$, $t_2=t_3$, and $\neg t_1=t_3$ are all in $\Sigma$. But then
$\Sigma$ is not finitely satisfiable.

\smallskip
So $\sim$ is an equivalence relation.  Next we show that it is a congruence
relation.

Suppose $t_1,\ldots,t_n,t_1^\prime,\ldots,t_n^\prime$ are variable free
terms and $f$ is an n-ary operation symbol.

If $t_1\sim t_1^\prime,\ldots,t_n\sim t_n^\prime$ then
$f(t_1,\ldots,t_n)\sim f(t_1^\prime,\ldots,t_n^\prime)$.

This amounts to saying if
$$\{t_1=t_1^{\prime},\ldots,t_n=t_n^{\prime}\}\subseteq\Sigma,$$
then
$$f(t_1,\ldots,t_n)= f(t_1^\prime,\ldots,t_n^\prime)\in\Sigma.$$
But again since $\Sigma$ is complete we would have
$$f(t_1,\ldots,t_n)\not= f(t_1^\prime,\ldots,t_n^\prime)\in\Sigma$$
but
$$\{t_1=t_1^{\prime},\ldots,t_n=t_n^{\prime},
f(t_1,\ldots,t_n)\not= f(t_1^\prime,\ldots,t_n^\prime)\}$$
has no models and
so $\Sigma$ wouldn't be finitely satisfiable.

By a similar argument: Suppose
$t_1,\ldots,t_n,t_1^\prime,\ldots,t_n^\prime$ are variable free
terms and $R$ is an n-ary operation symbol.

If $t_1\sim t_1^\prime,\ldots,t_n\sim t_n^\prime$ then
$R(t_1,\ldots,t_n)\in\Sigma$ iff
$R(t_1^\prime,\ldots,t_n^\prime)\in\Sigma$.

This shows the canonical model is well defined.

\medskip
Now we prove
by induction on the number of logical symbols that for any $L$ sentence
$\theta$,
$$\moda\models \theta \iff \theta \in \Sigma.$$

The atomic formula case is by definition.

\smallskip
$\neg$:
$$\moda\models\neg\theta$$
\centerline{iff}
$$\mbox{ not }\moda\models\theta$$
\centerline{iff(by induction)}
$$\mbox{ not }\theta\in\Sigma$$
\centerline{iff(by completeness)}
$$\neg\theta\in\Sigma.$$

\smallskip
$\disj$:
$$\moda\models(\theta\disj\rho)$$
\centerline{iff}
$$\moda\models\theta \rmor \moda\models\rho$$
\centerline{iff(by induction)}
$$\theta\in\Sigma \rmor \rho\in\Sigma$$
\centerline{iff}
$$(\theta\disj\rho)\in\Sigma.$$
This last ``iff'' uses completeness and finite satisfiability of $\Sigma$.
If $(\theta\disj\rho)\notin\Sigma$ then by
completeness $\neg(\theta\disj\rho)\in\Sigma$ but
$\{\theta,\rho,\neg(\theta\disj\rho)\}$ has no model.
Conversely if $\theta\notin\Sigma$ and $\rho\notin\Sigma$,
then by completeness $\neg\theta\in\Sigma$ and $\neg\rho\in\Sigma$, but
$\{(\theta\disj\rho),\neg\theta,\neg\rho\}$ has no model.

\smallskip
$\exists$:

$\moda\models\exists x\theta(x)$
implies
there exists $a\in A$ such that $\moda\models\theta(a)$.
This implies (by induction) $\theta(a)\in\Sigma$.
Which in turn implies  $\exists x\theta(x)\in\Sigma$,
since otherwise $\neg\exists x\theta(x)\in\Sigma$ but
$\{\neg\exists x\theta(x),\theta(a)\}$ has no model.
Hence $\moda\models\exists x\theta(x)$ implies $\exists x\theta(x)\in\Sigma$.

For the other direction suppose $\exists x\theta(x)\in\Sigma$.  Then
since $\Sigma$ is Henkin for some constant symbol $c$
we have $(\exists x\theta(x))\implies\theta(c)\in\Sigma$.  Using completeness
and finite satisfiability we must have $\theta(c)\in\Sigma$.
By induction $\moda\models\theta(c)$ hence $\moda\models\exists x\theta(x)$.
Hence $\exists x\theta(x)\in\Sigma$ implies
$\moda\models\exists x\theta(x)$.
\qed

\par\bigskip
\noindent {\bf Compactness Theorem.}  For any language $L$ and
set of $L$ sentences $\Sigma$, if every finite subset of
$\Sigma$ has a model, then $\Sigma$ has a model.

\bigskip

\proof
First Henkinize $\Sigma$, then complete it.  Take its canonical model.
Then reduct it back to an $L$-structure.
\qed

\bigskip

\begin{center}
 Problems
\end{center}

\prob  We say that $T$ a set of $L$ sentences is an \dex{$L$ theory}
iff there exists a set $\Sigma$ of $L$ sentences such
that $T$ is the set of all $L$ sentences true in every model
of $\Sigma$.  In this case we say that $\Sigma$ is an
\dex{axiomatization} of $T$ or that
$\Sigma$ \dex{axiomatizes} the theory $T$.
Prove that any theory axiomatizes itself.
Give an axiomatization of the theory of partially ordered sets.
The theory of groups is just the set of all sentences of
group theory which are true in every group.  Give an axiomatization of
the theory of abelian groups.

\prob  Suppose that $T$ is a theory with a model.
Show that $T$ is complete iff
for every sentence $\theta$
in the language of $T$ either  every model of $T$ is a model
of $\theta$ or every model of $T$ is a model
of $\neg\theta$.

\prob A theory $T$ is \dex{finitely axiomatizable} iff
there is a finite $\Sigma$
which axiomatizes $T$.  Let $T$ be the set of sentences in the language of
pure equality which are true in every infinite structure. Prove that
$T$ is a theory by finding axioms for it.  Show that no finite set
of these axioms axiomatizes $T$.

\prob
The theory of $\moda$, ${Th(\moda)}$,  is defined as follows:
\sdex{$Th(?moda)$}
$$Th(\moda)=\{\theta:\theta \mbox{ is an } L \mbox{ sentence and }
\moda\models\theta\}$$
Prove that $Th(\moda)$ is a complete theory.

\prob   $\moda$ and $\modb$ are \dex{elementarily equivalent} (written
${\moda\equiv\modb}$) iff $Th(\moda)=Th(\modb)$.
\sdex{$?moda?equiv?modb$}
Show that $\moda\equiv\modb$ iff for every sentence
$\theta$ if $\moda\models\theta$, then $\modb\models\theta$.

\prob
Suppose $T$ is a theory with a model.
Show the following are equivalent:
\begin{itemize}
\item $T$ is complete
\item any two models of $T$ are elementarily
equivalent
\item $T=Th(\moda)$ for some model of $T$
\item $T=Th(\moda)$ for all models of $T$
\end{itemize}

\prob Show that for any set of sentences  $\Sigma$ and sentence $\theta$,
every model of $\Sigma$ is a model of $\theta$ iff there exists a finite
$\Sigma^\prime\subseteq\Sigma$
such that every model of $\Sigma^{\prime}$ is a model of $\theta$.

\prob Suppose that $T$ is a finitely axiomatizable theory and
$\Sigma$ is any axiomatization of $T$.
 Show that some finite $\Sigma^\prime \subseteq \Sigma$
 axiomatizes $T$.

\prob (Separation) Let ${{\cal M}(T)}$
be the  \sdex{${?cal M}(T)$}
class of all models of $T$.
Suppose $T$ and $T^\prime$  are theories in a language $L$ and
 ${\cal M}(T) \cap {\cal M}(T^\prime ) = \emptyset $.  Show that
 for some $L$ sentence
 $\theta$,
 $${\cal M}(T) \subseteq {\cal M}(\theta) \rmand
  {\cal M}(\theta) \cap {\cal M}(T^\prime ) = \emptyset.$$

\prob Let $L$ be a first order language and suppose $T_i, i\in I$
are theories in $L$ such that every $L$ structure is a model of exactly
one of the $T_i$'s.  If $I$ is finite does it
then follow
that each of the $T_i$'s is finitely axiomatizable? What about infinite
$I$?

\prob Show that any theory with arbitrarily large finite models
must have an infinite model.

\prob Suppose that $T$ is an $L$-theory with an infinite model.
Let $X$ be any nonempty set and $\{c_x:x\in X\}$ be new constant
symbols not appearing in $L$.  Prove that the set of
sentences $T\cup\{c_x\not=c_y:x,y\in X,x\not=y\}$ has a model.
Prove that any theory with an infinite model has an uncountable
model.

\prob
Let the
language of fields be $\{ +, \cdot, 0 , 1 \}$.
Give an axiomatization for fields of characteristic $0$.
Show that this theory is not
finitely axiomatizable.  Show that for any sentence $\theta$
true in all fields of
characteristic $0$ there is a $k$ such that $\theta$ is true in all
fields of characteristic $ > k$.

\prob  The cardinality of a model is the number of elements in its
universe.
Give an example of a theory $T$  such that
for all $n$,
\par\centerline{$T$ has a model of cardinality $n$ iff $n$ is even.}
\par\noindent Can you find a finitely axiomatizable $T$?

\prob Is there a finitely axiomatizable theory with only infinite models?
What if the language consists of a single unary operation symbol? What
about the languages which contain only unary relation symbols?

\prob Suppose that $T$ is any theory in a language which includes a binary
relation symbol $\leq$ such that for every model $\moda$ of $T$
 $\leq^{\moda}$ is a
linear order. Show that if $T$ has an infinite model then $T$ has a model
$\modb$
such that there is an order embedding of the rationals into
 $\leq^{\modb}$, i.e., there is a function $f:{\Bbb Q}\to B$ such
 that $p\leq q$ iff $f(p)\leq^{\modb} f(q)$ for every $p,q\in {\Bbb Q}$.

\prob Let $\moda =( A,\eq )$ be the equivalence relation with
exactly one equivalence class of cardinality $n$ for each $n=1,2,\ldots$
and no infinite equivalence classes.
Show that there exists $\modb \equiv \moda$ which has infinitely
many infinite equivalence classes.

\prob  A theory is \dex{consistent} iff it has a model, i.e., it is realizable.
Let $T$ be a finitely axiomatizable theory with only a countable
number of complete consistent extensions (ie $T^\prime\supseteq T$)
in the language
of $T$.  Prove that one of these complete
consistent extensions is finitely axiomatizable.

\prob For each of the following prove or give a counterexample:
\begin{enumerate}
\item Let $T_n$ for $n\in\omega$ be finite sets of sentences and $S$ be
a finite set of sentences.
 Assume for all $n\in\omega$ that
$ T_n\subseteq T_{n+1}$, and there is a model of $T_n$ which is not a model
of $S$.  Then there is a model of $\cup_{n\in\omega}T_n $ which
is not a model of $S$.

\item Let $S_n$ and $T_n$ for $n\in\omega$
be finite sets of sentences.  Assume that for all
$n\in\omega,S_n \subseteq S_{n+1},
T_n\subseteq T_{n+1}$, and there is a model of $T_n$ which is not a model
of $S_n$.  Then there is a model of $\cup_{n\in\omega}T_n $ which
is not a model of $\cup_{n\in\omega}S_n $.
\end{enumerate}

\prob Let $ T_0 \subseteq T_1 \subseteq T_2 \subseteq \ldots $
be a sequence of
$L$ theories such that for each $ n \in \omega $ there exists a model of
$T_{n}$ which is not a model of $T_{n+1}$.
 Prove that $\cup_{n\in\omega}T_n$
 is not finitely axiomatizable. If $L$ is finite, prove that
 $\cup_{n\in\omega}T_n$ has an infinite model.

\sector{Lowenheim-Skolem Theorems}

The first version of the Lowenheim-Skolem Theorem was proved in 1915.
The final version that is presented here was developed by Tarski
in the 1950's.

\begin{lemma} The number of $L$ formulas is $|L|+\aleph_0$. \label{low1}
\end{lemma}
\proof
There are only countably many logical symbols.  Hence if
$\kappa=|L|+\aleph_0$ then every formula may be regarded as
an element of $\kappa^{<\omega}$ and we know
$|\kappa^{<\omega}|=\kappa$.
\qed

\bigskip
\par\noindent
{\bf Theorem}
{\em
For any theory $T$ in a language $L$ if $T$ has
a model, then it has one of cardinality less than or
equal to $|L|+\aleph_0$.
}
\proof
The canonical model of the completion of
Henkinization of $T$ has cardinality $\leq|L|+\aleph_0$.
It is enough to see that the language of the Henkinization
of $T$ has cardinality $\leq|L|+\aleph_0$, since the
canonical models universe is the set of equivalence classes
of variable free terms.
The Henkinization language is $\cup_{n<\omega}L_n$ where
$L_0=L$ and $L_{n+1}$ is $L_n$ plus one more constant symbol
for each formula of $L_n$.  But by Lemma \ref{low1}
there are $|L_n|+\aleph_0$ formulas of $L_n$.  It follows
that if $\kappa=|L|+\aleph_0$, then $|L_n|=\kappa$ for
all $n$ and so $|\cup_{n<\omega}L_n|=\kappa$.
\qed

${{\moda} \substruc {\modb}}$
\sdex{${?moda} ?substruc {?modb}$}
 means that ${\moda}$ is a
\dex{substructure}
 of  ${\modb}$;    equivalently
${\modb}$ is an \dex{extension} or \dex{superstructure}
 of ${\moda}$.  This means that
both structures are in the same language $L$, $A  \subseteq  B$, for
 every   $n$-ary
relation symbol $R$ of $L$, $ R^{\moda} = R^{\modb} \cap A^n$,
for every $n$-ary function symbol $f$ of $L$,
$f^{\moda} = f^{\modb} \res A^n $,
and for every constant symbol $c$ of $L$,
$ c^{\moda} = c^{\modb} $.

${{\moda} \elemsub {\modb}}$
\sdex{${?moda} ?elemsub {?modb}$}
means that ${\moda}$ is an
\dex{elementary substructure}
 of ${\modb}$;   equivalently ${\modb}$  is an \dex{elementary extension}
 of ${\moda}$.  This means that ${\moda} \subseteq {\modb}$
and for every formula $\theta(x_1,x_2,\ldots,x_n)$
 of the language $L$ and for every $a_1,a_2,\ldots,a_n \in A$
we have
$$({\moda },a)_{a\in A} \models \theta (c_{a_1},c_{a_2},\ldots,c_{a_n})
      \mbox{ iff }
  ({\modb },a)_{a\in A} \models \theta (c_{a_1},c_{a_2},\ldots,c_{a_n})$$
To ease the notational complexity we will write
${\moda\models\theta(a_1,\ldots,a_n)}$
\sdex{$?moda?models?theta(a_1,?ldots,a_n)$}
 instead of
$({\moda },a)_{a\in A} \models \theta (c_{a_1},c_{a_2},\ldots,c_{a_n})$.
It should be kept in mind that the language $L$ may have no constant
symbols in it.

Example.  Let $\moda=(\omega,<)$ and let $\modb=(Evens,<)$.  Then
$\modb\subseteq\moda$ but $\modb$ is not an elementary substructure
of $\moda$.  This is because $\moda\models\exists x\; 0<x<2$ but
$\modb\models\neg\exists x\; 0<x<2$.

\begin{lemma} (Tarski-Vaught criterion) Suppose $\moda\subseteq\modb$
are $L$ structures and
for every $L$ formula $\theta(x,y_1,y_2,\ldots,y_n)$, for all
$a_1,a_2,\ldots,a_n \in A$, and $b \in B$
$$ \modb\models\theta(b,a_1,a_2,\ldots,a_n)\mbox{ implies there exists }
 a \in A\,\,\,
\modb\models\theta(a,a_1,a_2,\ldots,a_n).$$
Then $\moda\elemsub\modb$.
\end{lemma}
\proof
The proof is by induction on the number of logical symbols in the
formula $\theta$.   The atomic formula case is trivial because
$\moda$ is a substructure of $\modb$.

$\neg$: $\moda\models\neg\theta$ iff not $\moda\models\theta$
iff(by induction) not $\modb\models\theta$ iff $\modb\models\neg\theta$.

$\disj$: $\moda\models(\theta\disj\rho)$ iff $\moda\models\theta$ or
$\moda\models\rho$
iff(by induction)  $\modb\models\theta$ or $\modb\models\rho$
iff $\modb\models(\theta\disj\rho)$.

$\exists$: $\moda\models\exists x \theta(x,a_1,\ldots,a_n)$ implies
there exists $a\in A$ such that $\moda\models \theta(a,a_1,\ldots,a_n)$
which implies (by induction) $\modb\models \theta(a,a_1,\ldots,a_n)$.
For the other direction we use the criterion.
$\modb\models\exists x \theta(x,a_1,\ldots,a_n)$ implies
there exists $a\in A$ such that
$\modb\models\theta(a,a_1,\ldots,a_n)$. Hence by induction
$\moda\models\theta(a,a_1,\ldots,a_n)$ and so
$\moda\models\exists x \theta(x,a_1,\ldots,a_n)$.
\qed

\begin{lemma} Suppose $X\subseteq B$ and $|X|=\kappa \geq |L|+\aleph_0$
where $\modb$ is an $L$ structure. Then there exists
$X^*\supseteq X, |X^*|=\kappa$,
and
for every formula $\theta(x,y_1,y_2,\ldots,y_n)$, for all
$a_1,a_2,\ldots,a_n \in X$, and $b \in B$
$$ \modb\models\theta(b,a_1,a_2,\ldots,a_n)\mbox{ implies }
\exists a \in X^*\,\,
\modb\models\theta(a,a_1,a_2,\ldots,a_n)$$
\end{lemma}
\proof
Fix $\leq$ a wellordering of $B$.
For any $L$-formula $\theta(x,y_1,\ldots,y_n)$ and
$$a_1,\ldots,a_n\in B,$$ define $a_{\theta(x,a_1,\ldots,a_n)}\in B$
to be the $\leq$ least element $b$ of $B$ such that
$\modb\models\theta(b,a_1,\ldots,a_n)$ if one exists otherwise
let it be arbitrary.
Now let $X_0=X$, $L_0=L$, and for any $m<\omega$ let
$X_{m+1}=$
$$\{a_{\theta(x,a_1,\ldots,a_n)}:\theta(x,y_1,\ldots,y_n)
\mbox{ is an  } L_m\mbox{ formula and }\{ a_1,\ldots,a_n\}\subseteq X_m\}$$
and let $L_{m+1}\supseteq L_m$ be the language with all these new
constant symbols adjoined.  Clearly if
$X_m$ and $L_m$ have cardinality $\kappa$ then so do
$X_{m+1}$ and $L_{m+1}$, since $|\kappa^{<\omega}|=\kappa$.
Let $X^*=\cup_{m<\omega}X_m$, then it has cardinality $\kappa$ since
it is the countable union of sets of cardinality $\kappa$.
For every  formula $\theta(x,a_1,a_2,\ldots,a_n)$ there exist
some $m<\omega$ such that $\{a_1,a_2,\ldots,a_n\}\subseteq X_m$ and
so the criterion for $\theta(x,a_1,a_2,\ldots,a_n)$ is satisfied
at stage $m+1$.
\qed

\bigskip

Definition: ${|{\moda}|}$
\sdex{$|{?moda}|$}
is the cardinality of
 the universe $A$ of ${\moda}$.

\par\bigskip
\noindent {\bf Downward Lowenheim-Skolem Theorem.}
Suppose ${\modb}$  is an infinite structure
in the language $L$, $\kappa$ is a cardinal  such that
$\aleph _0 + |L|  \leq  \kappa \leq |{\modb}|$, and
$X\subseteq B$ such that $|X|\leq \kappa$.  Then
there is a structure ${\moda}$
such that  ${\moda} \elemsub {\modb}$, $X\subseteq A$,
and $|{\moda}|=\kappa$.

\bigskip

\proof
By the lemma there exists $X^*\supseteq X$ of cardinality $\kappa$ satisfying
the criterion.  But note that the criterion implies that
$X^*$ is closed under the operations of $\modb$.  (Just look
at the sentence $\exists x\;\; x=f(a_1,\ldots,a_n)$.)  Consequently there
is a substructure $\moda$ of $\modb$ with universe $A=X^*$.
By the Tarski-Vaught criterion $\moda\elemsub\modb$.
\qed

\par\bigskip

${{\moda} \isom {\modb}}$
\sdex{${?moda} ?isom {?modb}$}
 means that  ${\moda}$ and ${\modb}$ are
\dex{isomorphic}, that is, there is a bijection $ j : A \to B$ such that
for every n-ary relation symbol R and for every $a_1,a_2,\ldots,a_n \in A$,
$$\langle a_1,a_2,\ldots,a_n\rangle \in R^{\moda}
\mbox{ iff }
\langle j(a_1),j(a_2),
\ldots,j(a_n)\rangle \in R^{\modb}
$$
and
for every n-ary function symbol $f$ and for every $a_1,a_2,\ldots,a_n \in A$,
$$j(f^{\moda}(a_1,a_2,\ldots,a_n))=
f^{\modb}(j(a_1),j(a_2),\ldots,j(a_n))$$
(for n=0 this means that for every constant symbol $c$,
$j( c^{\moda}) = c^{\modb}$.)

\begin{lemma}
Suppose $j$ is an isomorphism between the $L$ structures
$\moda$ and $\modb$.  Then for any $L$ formula
$\theta(x_1,x_2,\ldots,x_n)$ and any
$a_1,a_2,\ldots,a_n \in A$,
$$\moda\models \theta(a_1,a_2,\ldots,a_n)
\mbox{ iff }\modb\models \theta(j(a_1),j(a_2),\ldots,j(a_n))$$
\end{lemma}
\proof
First show by induction that for every $L$-term
$\tau(x_1,\ldots,x_n)$ and sequence $a_1,\ldots,a_n\in A$
that
$$j(\tau^\moda(a_1,\ldots,a_n))=\tau^\modb(j(a_1),\ldots,j(a_n)).$$
The proof of the lemma is by induction on the logical complexity of $\theta$.
For atomic formula it follows from the definition
of isomorphism.  The propositional steps are easy and
the quantifier step is handled by using that $j$ is onto.
\qed

A map $j:\moda\to\modb$ is an
\dex{elementary embedding} iff it is an
isomorphism of $\moda$ with an elementary substructure of $\modb$.
We write ${\moda\elemsub^j\modb}$.
\sdex{$?moda?elemsub^j?modb$}

\begin{lemma}
 A map  $j:\moda\to\modb$ is an elementary embedding iff
for any formula $\theta(x_1,x_2,\ldots,x_n)$ and any
$a_1,a_2,\ldots,a_n \in A$,
$$\moda\models \theta(a_1,a_2,\ldots,a_n)
\mbox{ iff }\modb\models \theta(j(a_1),j(a_2),\ldots,j(a_n))$$
\end{lemma}
\proof
Just unravel the definitions.
\unitlength=1.00mm
\linethickness{0.4pt}
\begin{picture}(55,20)
\put(0,0){\framebox(55,20)[]{}}
\put(10.00,10.00){\circle{10.00}}
\put(38.00,10.00){\circle{10.00}}
\put(39.00,10.00){\circle{15.00}}
\put(11.00,10.00){\vector(1,0){26.00}}
\put(24.00,13.00){\makebox(0,0)[cc]{$j$}}
\put(10.00,10.00){\makebox(0,0)[cc]{$\moda$}}
\put(43.00,13.00){\makebox(0,0)[cc]{$\modb$}}
\put(39.00,10.00){\makebox(0,0)[cc]{$\moda^\prime$}}
\end{picture}
\qed

\begin{lemma}
Suppose that $\moda\elemsub^j\modb$, then there exists a structure
$\modb^\prime$ isomorphic to $\modb$
such that $\moda\elemsub\modb^{\prime}$.
Furthermore $j$ is the restriction of this isomorphism to $A$.
\end{lemma}
\proof
Let $B^\prime$ be a superset of $A$ such that the map
$j$ can be extended to a bijection $j:B^{\prime}\to B$, (which
we also will call $j$).  Now define $f^{\modb^{\prime}}$ and
$R^{\modb^{\prime}}$ in such away as to make $j$ an isomorphism.
This means that
$$R^{\modb^{\prime}}=\{(b_1,\ldots,b_n):
                      (j(b_1),\ldots,j(b_n))\in R^{\modb}\}$$
and
$$f^{\modb^{\prime}}(b_1,\ldots,b_n)=
           j^{-1}(f^{\modb}(j(b_1),\ldots,j(b_n))$$
for all $b_1,\ldots,b_n\in B^\prime$.
Now check that $\modb^{\prime}$ works.
\qed

The \dex{elementary diagram} of $\moda$ is defined as follows:
$$D(\moda)=\{\theta:\theta \mbox{ is an } L_{\moda} \mbox{ sentence and }
(\moda,a)_{a\in A}\models\theta\}$$
This means that ${D(\moda)}$
\sdex{$D(?moda)$}
is the theory of $\moda$ with constants
adjoined for each element of the universe.

\begin{lemma} If $\moda$ is an $L$ structure and  $\modb$ is an $L_{\moda}$
structure such that
$\modb\models D(\moda)$, then
there is a $j$ such that $\moda\elemsub^j\modb\res  L$.
\end{lemma}
\proof
Define $j:A\to B$ by $j(a)=c_a^{\modb}$.
\qed

\par\bigskip
\noindent {\bf Upward Lowenheim-Skolem Theorem.}  For any infinite structure
${\moda}$ in the language $L$ and cardinal $\kappa$ such that
$|{\moda}| + |L| \leq \kappa$, there is a structure ${\modb}$
such that  ${\moda} \elemsub {\modb}$ and $|{\modb}|$=$\kappa$.

\bigskip

\proof
Let $\Sigma=D(\moda)\cup\{c_{\alpha}\not=c_{\beta}:\alpha,\beta\in\kappa,
\alpha\not=\beta\}$ where the $c_{\alpha}$ are new constant
symbols.

$\Sigma$ is finitely satisfiable.  To see this
let  $F\subseteq\Sigma$ be  finite.  Then there exists a
finite $G\subseteq\kappa$
such that $F\subseteq D(\moda)\cup
\{c_{\alpha}\not=c_{\beta}:\alpha,\beta\in G,
\alpha\not=\beta\}$.  Since the model $\moda$  is infinite we
can choose distinct elements of $a_\alpha\in A$ for $\alpha\in G$ and
then  $(\moda,a_{\alpha})_{\alpha\in G}\models F$.

It follows from the compactness theorem that $\Sigma$ has a model.
Since the language of $\Sigma$ has cardinality $\kappa$ by the
downward Lowenheim Skolem theorem $\Sigma$ has a model
$\modc$ of cardinality $\kappa$.  By the lemma there exists $j$ such that
$\moda\elemsub^j\modc\res L$.  By the other lemma $\modc\res L$ is
isomorphic to a model $\modb$ such that $\moda\elemsub\modb$.
\qed

\begin{center}
  Problems
\end{center}

\prob If an $L$ theory $T$ has an infinite model, then for any
$\kappa \geq
\aleph_0 +|L|$,  $T$ has a model of cardinality $\kappa$.

\prob Let $T$ be a consistent theory in a countable language.  Suppose
that $T$ has only countably many countable models up to isomorphism.
Show that for some sentence $\theta$ in the language of $T$ that
$T\cup\{\theta\}$ axiomatizes a complete consistent theory.

\prob Show that if $\moda \isom \modb$, then $\moda\equiv\modb$.
Is the converse true?  Show that if $\moda\elemsub\modb$,
then $\moda\equiv\modb$.  Is the converse true?

\prob Assume that $\moda\subseteq \modb$.
Show that $\moda\elemsub\modb$ iff for any $n$ and
$a_1,a_2,\ldots,a_n\in A$
we have $$(\moda,a_1,a_2,\ldots,a_n)\equiv (\modb,a_1,a_2,\ldots,a_n).$$

\prob A subset $X$ of a model $\moda$ is \dex{definable} iff for some
formula $\theta(x)$ in the language of $\moda$
$$ X = \{ b \in A : \moda \models \theta(b)\}.$$
Show that if $U\subseteq A$ is definable then for any isomorphism
$j:\moda\to\moda$ we have $U=\{j(x):x\in U\}$.
Show that the only definable subsets
of $(\Bbb Z,\leq)$ are the empty set and $\Bbb Z$ the set of all integers.

\prob A theory $T$ is \dex{categorical} in power $\kappa$ iff every two
models of $T$ of cardinality $\kappa$ are isomorphic.
 Suppose that $T$ is a consistent theory that has
only infinite models and for some infinite cardinal $\kappa \geq |L|$,
$T$ is $\kappa$ categorical.  Show that $T$ must be
complete. This is called the \dex{{?L}os-Vaught test}.
Is the assumption of no finite models  necessary?

\prob Let the language of $T_s$ be $\{S,c\}$ and
\begin{enumerate}
\begin{enumerate}
\item $ \theta   \equiv  \forall x \forall y \, (S(x)=S(y) \implies x=y) $
\item $ \rho     \equiv  \forall x \, (x \not= c \iff \exists y  S(y)=x) $
\item $ \psi _n  \equiv  \forall x \, S^n(x) \not= x$, where
$S^n(x)$ abbreviates $S(S(\cdots S(x)\cdots))$ where we have $n$ $S$'s.
\item  Let $T_s$ be the theory axiomatized by
$\{ \theta, \rho \} \cup \{ \psi_n : n = 1,2,\ldots \}$.
\end{enumerate}
\end{enumerate}
Show that
$\modn_{Sc} = (\omega,Sc,0)$ is a model of $T_s$,
where $Sc$ is the successor function ($Sc(n)=n+1$). \sdex{$?modn_{Sc}$}
Show that $T_s$ is not finitely axiomatizable.

\prob Is $T_s$ categorical in some power?  Show that $T_s$ is complete.

\prob Let $\modn_{Sc}^*$ be the model
consisting of $(\omega,Sc,0)$ plus a disjoint copy of $({\Bbb Z},Sc)$.
Show that if a subset of
 $\modn_{Sc}^*$
is definable  then it must either contain all of the ${\Bbb Z}$ part
or none of it.
A set is \dex{cofinite} iff it is the complement of a finite set.
Show that a subset of  $\modn_{Sc}$ is definable
iff it is finite or cofinite.
\com{Use compactness to get a model where there exists an infinite
$n$ such that $\theta(n)\conj\neg\theta(S(n))$.}

\prob Show that $\leq$ is not definable in
$\modn_{Sc}$.

\prob Axiomatize the theory of algebraically closed fields of
characteristic $k$.  Show that this theory is categorical in every uncountable
power.
Hence this theory is complete.  Show that the reals are not a definable
subset of $(\Bbb C,+,\cdot,0,1)$ where $\Bbb C$ is the complex numbers.

\prob Let $\moda$ be a countable structure in a countable language
that includes a unary predicate symbol $U$ such that
$ U^{\moda}$  is  infinite. Show that there is a countable
 $\modb \elemsup \moda$ such that $ U^{\modb} \not=  U^{\moda}$.

\prob Axiomatize the theory of dense linear orders without endpoints
(DLO).  That is linear orders without greatest or least element such that
between any two distinct elements lies a third.  For example, the rationals
or the reals
under their usual order are dense linear orders.
(Cantor) Show that any two countable
dense linear orders are isomorphic.   Show that
 $(\Bbb Q,\leq)\elemsub(\Bbb R,\leq)$
where $\Bbb Q$ is the rationals and $\Bbb R$ is the reals.  Is this
true if we add to our structures $ + $ and $ \cdot $?

\prob Can a theory have exactly $\aleph_0$ nonisomorphic models?  Can
such a theory be in a finite language?
\com{There are uncountably many uncountable cardinals.}

\prob
Show that for any wellordering $\modb =(B,\leq)$ there is a
countable wellordering elementarily equivalent to it.

\prob Prove or disprove:
 $(\moda \elemsub \modc \conj \modb \elemsub
\modc \conj \moda \subseteq \modb) \implies
               \moda \elemsub \modb$.

\com{shuffle of the two point order or use DCLO
say Z, Z+Z and Z+2Z}

\prob Prove or disprove:
$(\moda \elemsub \modb\conj\moda \elemsub \modc\conj
               \modb \subseteq \modc) \implies
               \modb \elemsub \modc$.

\prob Prove or disprove:
$(\moda \subseteq \modb \conj \moda \isom \modb) \implies
               \moda \elemsub \modb $.

\prob Prove or disprove:
$(\moda \elemsub \modb \elemsub \modc) \implies \moda \elemsub \modc$.

\prob Prove or disprove: if $\moda$ is isomorphic to a substructure
of $\modb$ and $\modb$ is isomorphic to a substructure of $\moda$,
then $\moda\isom\modb$.

\prob Prove or disprove:
$(\exists j \;\;\moda\elemsub^j\modb\conj \exists k\;\;\modb\elemsub^k\moda)
\implies
\moda\isom\modb$.

\com{$\omega$ many $\omega_1$ sized equivalence relations, in other
model one more countable class}

\prob The \dex{atomic diagram} of an $L$ structure $\moda$ (written
$D_0(\moda)$) is the set of quantifier free sentences in
$D(\moda)$.   \sdex{$D_0(?moda)$}
Show that for any $L$ structure $\moda$ and $L_{\moda}$
structure $\modb$ that
$\modb\models D_0(\moda)$ iff $\moda$ is isomorphic to a substructure
of $\modb\res L$.

\prob Let $\moda$ and $\modb$ be two infinite linear orderings.  Show that
$\modb$ is isomorphic to a substructure of some elementary extension
of $\moda$.

\prob Let $\moda \elemsub \modb$ and
$\moda \subseteq \moda^\prime$ and $A^\prime\cap B=A$.
Show that there exists $\modb^\prime  \supseteq \modb$ such that
$\moda^\prime \elemsub \modb^\prime $.

\prob Find a structure $\moda$ in a countable language such that $\moda$ has
exactly $\omega_1$ elementary substructures and every proper elementary
substructure is countable.
\com{
$(\omega_1,\leq)$ is an uncountable
wellordering such that for every $\alpha\in\omega_1$ the set
$\{\beta:\beta<\alpha\}$ is countable. }

\prob Suppose that $\moda$ is a finite model and
 $\moda \equiv \modb$.  Show that
$\moda \isom \modb$.  Warning: the language of $\moda$ may be infinite.

\prob Suppose
$$ \moda_0 \elemsub \moda_1 \elemsub \moda_2 \elemsub \ldots $$
and let $\modb = \cup_{n < \omega} \moda_n$.  Show that for each
$n<\omega$, $\moda_n \elemsub \modb$. This is called Tarski's
\dex{Elementary Chain Lemma}.

\prob A family $F$ of structures is directed iff any two structures of
$F$ are elementary substructures of some other structure in $F$.
Show that every structure in a directed family is an elementary substructure
of the union.

\prob Suppose $T$ is a theory in a countable language which has
an infinite model.  Show that $T$ has a countable
model that is not finitely generated.  A structure $\moda$ is
\dex{finitely generated} iff there is a finite set $F\subseteq A$ such
that no proper substructure of $\moda$ contains $F$.

\prob (DCLO) Axiomatize the theory of discrete linear orders without
 endpoints.
That is linear orders without greatest or least elements and
any element has an immediate successor and predecessor, for example the
integers $\Bbb Z$ (negative and positive) under their usual ordering.
Let $\moda$ be any model of DCLO. Define for $x,y \in A \; x \eq y $ iff
there are only finitely many elements of $A$ between $x$ and $y$.  Show this
is an equivalence relation and each equivalence class is isomorphic
to $\Bbb Z$, call such a class a $\Bbb Z$ chain.
Describe any model of DCLO.
Show that for any countable model of DCLO and any two $\Bbb Z$ chains
in the model there is an a countable elementary extension with a
$\Bbb Z$ chain in between.  Show that any countable model of DCLO is an
elementary
substructure of a model with countably many  chains
ordered like the rationals, i.e. $\Bbb Z\times\Bbb Q$.
Show that DCLO is a complete theory.

\prob Let $\modn$ be the model $(\omega,+,\cdot,\leq,0,1)$.
If $\modn^*$
is any proper elementary extension of $\modn$ we refer to it as a
\dex{nonstandard model} of arithmetic.
The elements of $ N^* \setminus N$ are the infinite
integers of $\modn^*$.  Show that the \dex{twin prime conjecture} is true,
i.e. there are infinitely many $p$ such that $p$ and $p+2$  are prime iff
in every nonstandard model of arithmetic $\modn^*$ there is an infinite
$p$ such that $\modn^*  \models$ ``$p$ and $p+2$ are prime''.

\prob Let  $\modn^*$ be any countable nonstandard model of arithmetic.
Show that
$$( N^*,\leq^{*})\isom (\omega +(\Bbb Z \times \Bbb Q),\leq)$$

\prob Show that for any set $X$ of primes that there
exist a countable model $M\equiv\modn$ and an $x\in M$ such
that for every prime number $p$
$$p\in X \rmiff M\models \exists y \; x=py.$$
Show that $\modn$ has $\cc$ many pairwise nonisomorphic countable
elementary extensions.

\prob Let $\modn_{full}$ be the model with universe $\omega$ and having a
relation symbol $\underline {R}$ for every relation R
of any arity on $\omega$, operation symbol $\underline {f}$ for
every operation on $\omega$, and a constant symbol $\underline {n}$ for
each element of $\omega$.  This is known as \dex{full arithmetic}.
Show that the language of full arithmetic has cardinality $\cc$.  Show that
every proper elementary extension of full arithmetic has cardinality
$\geq \cc$.
\com{consider eventually different functions.}

\prob Let $g(n)$ be the number of distinct prime divisors or $n$, with
$g(0)=g(1)=0$. Prove that in every nonstandard model $\modn^*$ of full
arithmetic, there is an element $b$ such that in  $\modn^*$,
$$b>g(b)>g(g(b))>g(g(g(b)))>\cdots$$

\prob Find a set of sentences $T$ in an uncountable language such that
$T$ has arbitrarily large finite models but no countable infinite model.

\prob Find a theory $T$ in an uncountable language that has no
finite models and has exactly one countable model up to isomorphism, but is
not complete.

\com{let T have $\exists^{\geq n}$ and $c_\alpha=c_\beta \iff c_\delta=
c_\gamma$ for $\alpha\not=\beta$ and $\delta\not=\gamma$.}

\sector{Turing Machines and Recursive Functions}

A \dex{Turing machine} is a partial function $m$ such that for some
finite sets $A$ and
$S$ the domain of $m$ is a subset of $S \times A$ and range of $m$
is a subset  of $S \times A \times \{l,r\}$.
$$\partial m: S \times A \to S \times A \times \{l,r\}$$
We call $A$ the alphabet and $S$ the states.

For example, suppose  $S$ is the set of letters $\{a,b,c,\ldots,z\}$
and $A$ is the set of all integers less than seventeen, then
$$m(a,4)=(b,6,l)$$
would mean that when the machine $m$ is in state $a$
reading the
symbol $4$ it will go into state $b$, erase the symbol $4$ and
write the symbol $6$ on the tape square
where $4$ was, and then move left one square.

\bigskip

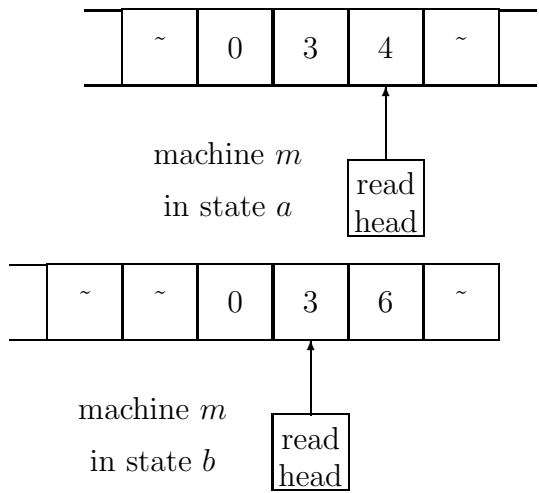
\unitlength=1.00mm
\special{em:linewidth 0.4pt}
\linethickness{0.4pt}
\begin{picture}(69.00,30.00)(-10,0)
\put(10.00,20.00){\framebox(10.00,10.00)[cc]{$\blank$}}
\put(20.00,20.00){\framebox(10.00,10.00)[cc]{$0$}}
\put(30.00,20.00){\framebox(10.00,10.00)[cc]{$3$}}
\put(40.00,20.00){\framebox(10.00,10.00)[cc]{$4$}}
\put(50.00,20.00){\framebox(10.00,10.00)[cc]{$\blank$}}
\put(40.00,0.00){\framebox(10.00,10.00)[cc]{}}
\put(45.00,2.00){\makebox(0,0)[cc]{head}}
\put(45.00,7.00){\makebox(0,0)[cc]{read}}
\put(45.00,10.00){\vector(0,1){10.00}}
\put(24.00,11.00){\makebox(0,0)[cc]{machine $m$}}
\put(24.00,4.00){\makebox(0,0)[cc]{in state $a$}}
\put(5.00,20.00){\line(1,0){5.00}}
\put(5.00,30.00){\line(1,0){5.00}}
\put(60.00,30.00){\line(1,0){5.00}}
\put(60.00,20.00){\line(1,0){5.00}}
\end{picture}

\bigskip

\unitlength=1.00mm
\special{em:linewidth 0.4pt}
\linethickness{0.4pt}
\begin{picture}(69.00,30.00)
\put(10.00,20.00){\framebox(10.00,10.00)[cc]{$\blank$}}
\put(20.00,20.00){\framebox(10.00,10.00)[cc]{$\blank$}}
\put(30.00,20.00){\framebox(10.00,10.00)[cc]{$0$}}
\put(40.00,20.00){\framebox(10.00,10.00)[cc]{$3$}}
\put(50.00,20.00){\framebox(10.00,10.00)[cc]{$6$}}
\put(60.00,20.00){\framebox(10.00,10.00)[cc]{$\blank$}}
\put(40.00,0.00){\framebox(10.00,10.00)[cc]{}}
\put(45.00,2.00){\makebox(0,0)[cc]{head}}
\put(45.00,7.00){\makebox(0,0)[cc]{read}}
\put(45.00,10.00){\vector(0,1){10.00}}
\put(24.00,11.00){\makebox(0,0)[cc]{machine $m$}}
\put(24.00,4.00){\makebox(0,0)[cc]{in state $b$}}
\put(5.00,20.00){\line(1,0){5.00}}
\put(5.00,30.00){\line(1,0){5.00}}
\end{picture}

\bigskip

If $(a,4)$ is not in the domain of $m$, then the machine halts.
This is the only
way of stopping a calculation.

Let $A^{<\omega}$ be the set of all finite strings from the alphabet A.
The Turing machine $m$ gives rise to a partial function
$M$ from $A^{<\omega}$ to $A^{<\omega}$ as follows.
$$\partial M:A^{<\omega}\to A^{<\omega}$$
We suppose that $A$ always contains
the blank
space symbol:
 $$\blank$$
and $S$ contains the starting state $a$.
Given any word $w$
from $A^{<\omega}$ we
imagine a tape with $w$ written on it and blank symbols everywhere else.  We
start the machine in state $a$ and reading the leftmost symbol of $w$.
A configuration consists of what is written on the tape, which square of
tape is being read, and the state the machine is in.
 Successive
configurations are obtained according to rules determined by $m$, namely
if the machine is in state $q$ reading symbol $s$ and
 $m(q,s)=(q^\prime,s^\prime,d)$ then
the next configuration has the same tape except the square we were reading
now has the symbol $s^\prime$ on it, the new state is $q^\prime$,
and the square being
read is one to the left if $d=l$ and one to the right if $d=r$.
If $(q,s)$ is not in the domain of $m$, then the computation halts and
$M(w)=v$ where $v$ is what is written on the tape when the machine halts.

Suppose $B$ is a finite alphabet that does not contain the blank
space symbol ($\blank$).   Then
a partial function
$$\partial f:B^{<\omega}\to B^{<\omega}$$
is a \dex{partial recursive function}
iff there
is a Turing machine $m$ with an alphabet $A\supseteq B$
such that $f = M\res B^{<\omega}$.  A partial recursive function
is \dex{recursive} iff it is total.  A function $f:\omega \to \omega$
is recursive iff it is recursive when considered as a map
from $B^{<\omega}$ to $B^{<\omega}$
where
$B=\{1\}$.  Words in $B$ can be regarded as numbers written
in base one, hence we identify the number $x$ with $x$ ones written
on the tape.

For example, the identity function is recursive, since it is computed by the
empty machine. The successor function is recursive since it is computed
by the machine:

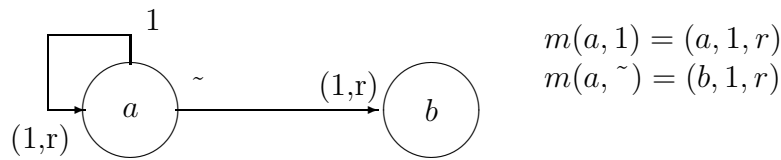
\unitlength=1.00mm
\special{em:linewidth 0.4pt}
\linethickness{0.4pt}
\begin{picture}(128.00,24.00)

\put(21.00,7.00){\vector(1,0){27.00}}
\put(15.00,17.00){\line(0,-1){4.00}}
\put(4.00,7.00){\vector(1,0){5.00}}
\put(4.00,7.00){\line(0,1){10.00}}
\put(4.00,17.00){\line(1,0){11.00}}

\put(18.00,19.00){\makebox(0,0)[cc]{1}}
\put(3.00,3.00){\makebox(0,0)[cc]{(1,r)}}
\put(24.00,10.00){\makebox(0,0)[cc]{$\;\tilde{}\;$}}
\put(44.00,10.00){\makebox(0,0)[cc]{(1,r)}}
\put(15.00,7.00){\circle{12.00}}
\put(55.00,7.00){\circle{12.00}}
\put(15.00,7.00){\makebox(0,0)[cc]{$a$}}
\put(55.00,7.00){\makebox(0,0)[cc]{$b$}}
\put(70,10){\shortstack[l]{ $m(a,1)=(a,1,r)$\\ $m(a,\blank)=(b,1,r)$}}
\end{picture}

In the diagram on the left, states  are represented by circles.
The arrows represent the \dex{state transition function} $m$.
For example, the horizontal arrow represents the fact that when
$m$ is in state $a$ and reads ($\blank$),  it
writes $1$, moves right, and goes into state $b$.

The set of strings of zeros and ones with
an even number of ones is recursive.
Its characteristic function (parity checker) can be computed
by the following machine:

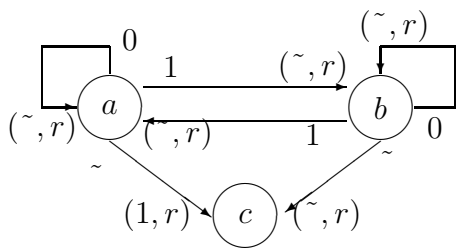
\unitlength=.90mm
\special{em:linewidth 0.4pt}
\linethickness{0.4pt}
\begin{picture}(140.00,40.00)
\put(15.00,025.00){\circle{10.00}}
\put(55.00,025.00){\circle{10.00}}
\put(35.00,009.00){\circle{10.00}}
\put(15.00,025.00){\makebox(0,0)[cc]{$a$}}
\put(55.00,025.00){\makebox(0,0)[cc]{$b$}}
\put(35.00,009.00){\makebox(0,0)[cc]{$c$}}
\put(20.00,028.00){\vector(1,0){30.00}}
\put(50.00,023.00){\vector(-1,0){30.00}}
\put(15.00,020.00){\vector(4,-3){15.00}}
\put(55.00,020.00){\vector(-4,-3){14.00}}
\put(15.00,030.00){\line(0,1){4.00}}
\put(15.00,034.00){\line(-1,0){10.00}}
\put(05.00,034.00){\line(0,-1){9.00}}
\put(05.00,025.00){\vector(1,0){5.00}}
\put(60.00,025.00){\line(1,0){6.00}}
\put(66.00,025.00){\line(0,1){9.00}}
\put(66.00,034.00){\line(-1,0){11.00}}
\put(55.00,034.00){\vector(0,-1){4.00}}
\put(18.00,035.00){\makebox(0,0)[cc]{$0$}}
\put(05.00,022.00){\makebox(0,0)[cc]{$(\blank,r)$}}
\put(24.00,031.00){\makebox(0,0)[cc]{$1$}}
\put(45.00,031.00){\makebox(0,0)[cc]{$(\blank,r)$}}
\put(45.00,021.00){\makebox(0,0)[cc]{$1$}}
\put(25.00,021.00){\makebox(0,0)[cc]{$(\blank,r)$}}
\put(13.00,015.00){\makebox(0,0)[cc]{$\blank$}}
\put(22.00,009.00){\makebox(0,0)[cc]{$(1,r)$}}
\put(47.00,009.00){\makebox(0,0)[cc]{$(\blank,r)$}}
\put(56.00,017.00){\makebox(0,0)[cc]{$\blank$}}
\put(63.00,022.00){\makebox(0,0)[cc]{$0$}}
\put(57.00,037.00){\makebox(0,0)[cc]{$(\blank,r)$}}
\put(80,10){\shortstack[l]{
$ m(a,0)=(a,\blank ,r)$       \\
$ m(a,1)=(b,\blank ,r)$       \\
$ m(b,0)=(b,\blank ,r)$       \\
$ m(b,1)=(a,\blank ,r)$       \\
$ m(a,\blank)=(c,1,r)$        \\
$ m(b,\blank)=(c,\blank,r)$   }}
\end{picture}

 The following problems are concerned with recursive functions and
predicates on $\omega$.

\bigskip

\prob Show that any constant function is recursive.

\prob  A binary function $f:\omega\times\omega\to\omega$ is recursive
iff there is a machine such that for any $x,y\in\omega$ if we input
$x$ ones and $y$ ones separated by ``,'', then the machine
eventually halts with $f(x,y)$ ones on the tape.
Show that $f(x,y) = x+y$ is recursive.

\prob Show that $g(x,y) = x y$ is recursive.

\prob Let $ x \dotminus y = max\{ 0,x-y \}.$ Show that $p(x)= x \dotminus 1$
is recursive.  Show that $ q(x,y)= x\dotminus y $ is recursive.

\prob Suppose $f(x)$ and $g(x)$ are recursive.  Show that $f(g(x))$
 is recursive.

\prob Formalize a notion of multitape Turing machine. Show that we get
the same set of recursive functions.

\prob Show that we get the same set of recursive functions even if we
restrict our notion of computation to allow only tapes that are infinite
in one direction.

\prob Show that the family of recursive functions is closed under arbitrary
compositions, for example $f(g(x,y),h(x,z),z)$.  More generally,
if $f(y_1,\ldots,y_m)$, $g_1(x_1,\ldots,x_n),\ldots$, and $g_m(x_1,\ldots,x_n)$
are all recursive, then so is
$$f(g_1(x_1,\ldots,x_n),\ldots,g_m(x_1,\ldots,x_n)).$$

\prob Define
$$ sgn(n) = \left\{
\begin{array}{ll}
 0 & \mbox{ if } n=0 \\
 1 & \mbox{ otherwise}
\end{array}
\right.$$
Show it is recursive.

\prob  A set is recursive iff its characteristic function is.
Show that the binary relation $x=y$ is recursive.  Show that the
binary relation $ x\leq y $ is recursive. \sdex{recursive set}

\prob Show that if $A \subseteq \omega$ is recursive then so is
$\omega\setminus A.$ Show that if $A$ and $B$ are recursive, then so are
$A\cap B$ and $A\cup B$.

\prob Suppose g(x) and h(x) are recursive and A is a recursive set. Show
that f is recursive where:
$$ f(x) =
\left\{
\begin{array}{ll}
 g(x) & \mbox{ if } x \in A \\
 h(x) & \mbox{ if } x \notin A
\end{array}
\right.$$

\prob Show that the set of even numbers is recursive.  Show that the
set of primes is recursive.

\prob Show that $e(x,y)=x^y$ is recursive.  Show that $f(x)=x!$ is recursive.

\prob Suppose that $h(z)$ and $g(x,y,z)$ are recursive. Define $f$
by recursion,
\begin{itemize}
\item $f(0,z)=h(z)$
\item $f(n+1,z) = g( n,z, f(n,z) )$.
\end{itemize}
Show that $f$ is recursive.

\prob We say that a set $ A \subseteq \omega$ is
\dex{recursively enumerable}
 (re) iff
it is the range of a total recursive function or the empty set.
Show that a set is re
iff it is the domain of a partial recursive function.

\prob  Show that every recursive set is re.  Show that a set is recursive
iff it and its complement are re.

\prob  Show that if $f$ is an increasing total recursive function then the
range of $f$ is recursive.

\prob Suppose that $f:\omega\to\omega$ and $g:\omega\to\omega$ are
recursive functions such that $f(m)<g(n)$ whenever $m<n$.
Prove that either the range of $f$ or the range of $g$ (or both)
is recursive.

\prob Let $f(n)$ be the $n^{th}$ digit after the ``.'' in the decimal
expansion of $e$. ($f(1)=7$, $f(2)=1$, $f(3)=8$).  Prove that the function
$f$ is computable.

\prob Show that there does not exist a total recursive function $f(n,m)$
 such
that for every total recursive function $g(m)$ there is an $n$ such that
$f(n,m)=g(m)$ for every $m$.

\bigskip

\begin{center}
 {\bf  Church's Thesis}
\end{center}

  One day in the 1930's Alonzo Church said ``Say, fellas, I
think that every function that is effectively computable is recursive.
I think we have captured this intuitive notion by this formal definition.
From now on why don't we just assume if we describe an effective
procedure it is possible to write done a Turing machine that does
it.  This will save a lot of time verifying silly details.''\footnote{
Actually these were probably not his exact words, in particular, he
was interested in some bizarre notion known as the lambda calculus.}
This philosophical position is known as  \dex{Church's Thesis}.

Here is an excerpt in support of Church's Thesis from
Alan M. Turing\footnote{
``On computable numbers, with an application to the Entscheidungsproblem'',
Proceedings of the London Mathematical Society,
2-32(1936), 230-265.}.  Note that Turing uses the word computer
for the person that is performing some effective procedure.

``Computing is normally done writing certain symbols on paper. We
   may suppose this is divided into squares like a child's arithmetic
   book.  In elementary arithmetic the two-dimensional character of the
   paper is sometimes used.  But such a use is always avoidable, and I
   think that it will be agreed that the two-dimensional character
   of paper is not essential for computation.  I assume then that the
   computation is carried out on one-dimensional paper, i.e., on a tape divided
   into squares.  I shall also suppose that the number of symbols which may be
   printed is finite.  If we were to allow an infinity of symbols, then
   there would be symbols differing to an arbitrarily small extent.
   The effect of this restriction of the number of symbols is not very
   serious. It is always possible to use sequences of symbols in the place
   of a single symbol.  Thus an Arabic numeral 17 or 9999999999999999999
   is normally treated as a single symbol.  Similarly in any European language
   words are treated as single symbols (Chinese, however, attempts to have
   an infinity of symbols).   The differences from our point of view between
   the single and compound symbols is that the compound symbols, if they
   are too lengthy, cannot be observed at one glance.  This is in accordance
   with experience.  We cannot tell at one glance whether
   9999999999999999999999999 and 99999999999999999999999999 are the same.

``The behavior of the computer at any moment is determined by the symbols
   which he is observing, and his `state of mind' at that moment.  We
   may suppose that there is a bound $B$ to the number of symbols or
   squares which the computer can observe at one moment.  If  he wishes
   to observe more, he must use successive observations.  We will also
   suppose that the number of states of mind which need be taken into
   account is finite.  The reasons for this are of the same character as
   those which restrict the number of symbols.  If we admitted an infinity
   of states of mind, some of them will be `arbitrarily close' and will
   be confused.   Again, the restriction is not one which seriously affects
   computation, since the use of more complicated states of mind can be
   avoided by writing more symbols on the tape.

``Let us imagine the operations performed by the computer to be split
   up into `simple operations' which are so elementary that it is not
   easy to imagine them further divided.   Every such operation consists
   of some change of the physical system consisting of the computer and
   his tape.  We know the state of the system if we know the sequence
   of symbols on the tape, which of these are observed by the computer
   (possibly with a special order), and the state of mind of the computer.
   We may suppose that in a simple operation not more than one symbol
   is altered.  Any other changes can be split up into simple changes of
   this kind.  The situation in regard to squares whose symbols may be
   altered in this way is the same as in regard to the observed squares.
   We may, therefore, without loss of generality, assume
   that the squares whose symbols are changed are always `observed'
   squares.

``Besides these changes of symbols, the simple operations must include
   changes of distribution of observed squares.  The new observed squares
   must be immediately recognizable by the computer.  I think it is
   reasonable to suppose that they can only be squares whose distance
   from the closest of the immediately previously observed squares
   does not exceed a certain fixed amount....

``The operation actually performed is determined, as has been suggested
   above, by the state of mind of the computer and the observed symbols.
   In particular, they determine the state of mind of the computer after
   the operation.''

\bigskip

Other evidence for Church's thesis is the fact that all
other ways people have come up with to formalize the notion of recursive
functions
(eg RAM machines, register machines, generalized recursive functions, etc)
can be shown to define the same set of functions.

In his paper Turing also proved the following remarkable theorem.

\bigskip
{\bf Universal Turing Machine Theorem}
{\em There is a partial  recursive function
$f(n,m)$ such that for every partial recursive function $g(m)$
there is an $n$ such that for every $m$,  $f(n,m)=g(m)$.}

\bigskip

Equality here means either both sides are
defined and equal or both sides are undefined.

\proof
Given the integer $n$ we first decode it as a sequence of integers
by taking its prime factorization,
$n=2^{k_1}3^{k_2}\cdots p_m^{k_m}$ ($p_m$ is the $m^{th}$ prime number).
Then we regard each integer $k_j$ as some character on the
typewriter (if $k_j$ too big we ignore it).
If the message coded by $n$ is a straight forward
description of a Turing machine, then we carry out the computation
this machine would do when presented with input $m$. If this
simulated computation halts with output $k$, then we halt with output $k$.
If it doesn't halt, then neither does our simulation.  If $n$ does not
in a straight forward way code the description of a Turing machine,
then we go into an infinite loop, i.e., we just never halt.
\qed

\bigskip

\prob Let $f$ be the universal function above and let
$K=\{ n : \langle n,n \rangle \in  dom(f)\}$. Show that $K$ \sdex{$K$}
 is re
but not recursive.

\prob The family of re sets can be uniformly enumerated by
$\langle W_e:e\in\omega\rangle$ where $W_e=\{n:(e,n)\in dom(f)\}$.
Show there exists a recursive function $d:\omega\to\omega$
such that for any $e\in\omega$ if $K\cap W_e=\emptyset$, then
$d(e)\notin K\cup W_e$.  This $d$ effectively witnesses that the
complement of $K$ is not re.

\com{Let $d(e)=e$.}

\sector{Completeness theorem}

The completeness theorem was proved by Kurt G\"{o}del in 1929.
To state the theorem we must  formally define the notion of proof.
This is not because it is good to give formal proofs, but rather so
that we can prove mathematical theorems about the concept of proof.

There are many systems for formalizing proofs in first order logic.
One of the oldest is the Hilbert-Ackermann style natural deduction system.
Natural systems seek to mimic commonly used
informal methods of proofs.   Another such system is the
 Beth semantic tableau method where a proof looks like
a tree of formulas.  This system is often used to teach undergraduates
how to do formal proofs.\footnote{Corazza, Keisler, Kunen, Millar,
 Miller, Robbin,
{\bf Mathematical Logic and Computability}.}

Other systems,  such as Gentzen style sequence calculus,
were invented with another purpose in mind, namely to analyze the
proof theoretic strength of various versions
of arithmetic.

Some  proof systems were invented to be easy
for a computer to use.  Robinson's resolution method is
popular with artificial intelligence people who try to get the
computer to prove mathematical theorems.   On the other hand it
is  hard for a human being to read a proof in this style.

\bigskip

The proof system we will use is constructed precisely to make the
completeness theorem easy to prove.

Definition of \dex{proof}:
If $\Sigma$ is a set of sentences and $\theta$ is any sentence, then
$\Sigma$ proves $\theta$ (and we write ${\Sigma\proves\theta}$)
\sdex{$?Sigma?proves?theta$}
iff there exists a finite sequence
$\theta_1,\ldots,\theta_n$ of sentences such that $\theta_n=\theta$ and
for each $i$ with $1\leq i\leq n$
\begin{enumerate}
\item $\theta_i$ is an element of $\Sigma$, or
\item $\theta_i$ a logical axiom, or
\item $\theta_i$ can be obtained by a logical rule from
the earlier $\theta_j$ for $j<i$.
\end{enumerate}

\medskip

Logical axioms and rules depend of course on the system being used.

\medskip

{\bf Completeness Theorem:} For any set of sentences
$\Sigma$ and sentence $\theta$
we have that:

\centerline{ $\Sigma\proves\theta$ iff every model of $\Sigma$
is a model of $\theta$.}

{\bf MM proof system:}
In the MM proof system (for Mickey Mouse) the logical axioms
and logical rules are the following.  We let all validities be
logical axioms, i.e.,  if a sentence $\theta$ is true in every
model, then $\theta$ is a logical axiom.  In the MM system there
is only one logical rule: \dex{Modus Ponens}.  This is the rule
that from $(\psi\implies\rho)$ and $\psi$ we can infer
$\rho$.  So in a proof if
some $\theta_i=(\psi\implies\rho)$ and some $\theta_j=\psi$
and $k>i,j$, then we can apply Modus Ponens to get $\theta_k=\rho$.

\bigskip
{\bf proof of completeness:}
The easy direction of the completeness theorem is often called
the \dex{soundness theorem}: if $\Sigma\proves\theta$ then
every model of $\Sigma$ is a model of $\theta$.  To check that MM
is sound is an easy induction on the length of the proof:
If $\theta_1,\ldots,\theta_n$ is a sequence of sentences
such that each $\theta_i$ is either an element of $\Sigma$, a logical
validity, or can be obtained by Modus Ponens from
the earlier $\theta_j$ for $j<i$; then by induction on
$j\leq n$ every model of $\Sigma$ is a model of $\theta_j$.

   Now we prove the other direction of the completeness theorem.
So suppose that every model of $\Sigma$ is a model of $\theta$, we
need to show that $\Sigma\proves\theta$. Since every model of
$\Sigma$ is a model of $\theta$ the set of sentences
$\Sigma\cup\{\neg\theta\}$ has
no models.  By the compactness theorem
there exists $\{\theta_1,\ldots,\theta_n\}\subseteq \Sigma$ such that
$\{\theta_1,\ldots,\theta_n,\neg\theta\}$ has no model.  This means
that
$$(\theta_1\implies(\theta_2\implies(\cdots\implies
(\theta_n\implies\theta)\cdots)))$$
is a validity and hence it is a logical axiom of system MM.  Now it is
easy to show that $\Sigma\proves\theta$: just write down
$\theta_1,\ldots,\theta_n$, and this logical axiom,
and then apply Modus Ponens $n$ times.
\qed

\bigskip
\noindent{\bf Corollary 1 of the Completeness Theorem:}

 If $\Sigma$ axiomatizes a theory $T$, then
$T=\{\theta:\Sigma\proves\theta\}$.

\bigskip

A set of sentences  $\Sigma$ is \dex{inconsistent} iff there exists
a propositional contradiction $\#$ such that
$\Sigma\proves\#$. All contradictions
are logically equivalent so we write them $\#$.
Note that if $\Sigma$ is inconsistent, then $\Sigma$ proves every
sentence, because if $\#$ is a contradiction, then $(\#\implies\theta)$
is a propositional validity, so if $\Sigma\proves\#$, then
$\Sigma\proves\theta$.
A set of sentences  $\Sigma$ is \dex{consistent} iff it is not inconsistent.

\bigskip
\noindent{\bf Corollary 2 of the Completeness Theorem:}
\par Every consistent set of sentences has a model.

\bigskip

The poor MM system went to the Wizard of OZ and said,
``I want to be more like all the other proof systems.''

And the Wizard
replied,
``You've got just about everything any other proof system has
and more.   The completeness theorem is  easy to prove in your system.
You have very few logical rules and logical axioms.   You lack
only one thing.  It is too hard for mere mortals
to gaze at a proof in your system and tell whether it really is
a proof.   The difficulty comes from taking all logical validities as
your logical axioms.''

The Wizard went on to give MM a subset $Val$ of logical validities
that is recursive and has the property that every
logical validity can be proved using only Modus Ponens from $Val$.

Let $L$ be the largest countable language you can think of, i.e.,
for each $n,m<\omega$
an operation symbol $f_{n,m}$ of arity $m$ and a relation symbol
$R_{n,m}$ of arity $m$. Note that we think of the $f_{n,0}$ as constant
symbols and the $R_{n,0}$ as propositional symbols.
The symbols of $L$ can be written
in a finite alphabet if we imagine that each $n<\omega$ is a string of
the symbols $\{0,1,2,\ldots,9\}$.
Since there are only finitely many logical symbols and variables
are written $x_n$ for $n<\omega$, the following definition makes sense.

A set $\Sigma$ of formulas in the language $L$ is recursive iff there
exists an effective procedure  which will decide whether
$\theta$ is in $\Sigma$.  That is, there exists a machine which when
given  any string $\theta$
of symbols in this finite alphabet will eventually halt and
say ``yes'' or ``no'' depending on whether $\theta$ is in $\Sigma$.
A set $\Sigma$ of formulas is recursively enumerable iff there exists
an effective procedure that will list all the formulas in $\Sigma$.

\prob  Show that the set of sentences of $L$ is recursive.

\com{
Take a course in how to write a compiler from the CS department.
}

\prob
Show that the set of propositional sentences of $L$ that are
validities is recursive.

\com{ Use truth tables.}

\bigskip

\prob An $L$ formula $\theta$ is in \dex{prenex normal form} iff
$$\theta=Q_1 x_1 \;Q_2x_2\cdots Q_n \; x_n\psi$$
where $\psi$ is quantifier free and each $Q_i$ is either $\forall$ or
$\exists$ (i.e.,  all quantifiers occur up front).
Show there exists a recursive map p from $L$-formulas to
prenex normal form $L$-formulas such that $\theta$ and $p(\theta)$ are
logically equivalent for every $\theta$.
Two formulas $\theta(x_1,\ldots,x_n),\psi(x_1,\ldots,x_n)$ are
logically equivalent iff
$\forall x_1\cdots\forall x_n(\theta(x_1,\ldots,x_n)\iff\psi(x_1,\ldots,x_n))$
is a logical validity.

\prob
Let $f$ is a new operation symbol not in appearing in $\theta(x,y)$.
Show $\forall x \exists y\theta(x,y)$ has a model iff
$\forall x\theta(x,f(x))$ has a model.

\prob  (Skolemization) An $L$ formula $\theta$ is
\dex{universal} iff it is in prenex normal form and all
the quantifiers are universal ($Q_i$ is $\forall$ for each $i$).
Show
there exists a recursive map s from $L$ sentences to universal
$L$ sentences such that
$$\theta \mbox{ has a model iff }s(\theta) \mbox{ has a model }$$
for every $\theta$ an $L$ sentence. \sdex{Skolemization}

\prob
Similarly an $L$ formula is \dex{existential} iff
it is in prenex normal form and all
the quantifiers are existential ( $Q_i$ is $\exists$ for each $i$).
Show there exists a recursive map $e$ from $L$ sentences to  existential
$L$ sentences such that

   \centerline{$\theta$ is a validity iff $e(\theta)$ is a validity}

\noindent for every $\theta$ an $L$ sentence.

\com{Consider  $\neg(s(\neg\theta))$.}

\prob  (Herbrand)
Suppose we have an existential $L$ sentence $\theta$ of the form
$$\exists x_1\ldots\exists x_n\;\psi(x_1,\ldots,x_n)$$
where
$\psi$ is quantifier free.
Show that $\theta$ is a validity iff for some
$m<\omega$ there exists a sequence
$$\vec{\tau_1},\vec{\tau_2},\ldots \vec{\tau_m}$$
such that each $\vec{\tau_i}$ is a $n$-tuple of variable
free terms of $L$ and
$$[\psi(\vec{\tau_1})\disj\psi(\vec{\tau_2})\disj\cdots\disj
\psi(\vec{\tau_m})] \mbox{ is a validity.}$$
\sdex{Herbrand}

\com{
One direction is easy.  For the other consider the smallest substructure
of a model of
$$\Sigma=\{\neg\psi(\vec{\tau}): \vec{\tau}
                 \mbox{ is an $n$-tuple of variable free terms}\}.$$
}

\prob Show that the set of quantifier free sentences of $L$ that
are validities is recursive.

\prob Prove that the set of logical validities of $L$ is
recursively enumerable.

\prob Find a recursive set of logical validities ${Val}$
\sdex{$Val$} in the
language $L$ such that every logical validity in $L$ can be proved
using only Modus Ponens and logical axioms from $Val$.

\com{If $\theta_0,\theta_1,\ldots,$ is a recursive enumeration of
sentences such that the length of $\theta_n$ is less than
$\theta_{n+1}$ for each $n$, then the set $\{\theta_0,\theta_1,\ldots,\}$
is recursive. }

\com{as well as all things of form
$(\theta_0\conj\cdots\conj\theta_n)\implies\theta_i)$}

\bigskip

{\bf In the following problems assume
all languages  are  recursive subsets of $L$.}

\bigskip

\prob  Show that for any recursive set $Q$ of sentences
the set of all $\theta$ such that $Q\proves \theta$ is
recursively enumerable.

\prob A theory is said to be \dex{recursively axiomatizable} iff it can
be axiomatized by a recursive set of sentences.
(Craig) Show that any recursively enumerable theory is
recursively axiomatizable.

\prob  A theory axiomatized by $\Sigma$ is \dex{decidable} iff the set
{$\{\theta:\Sigma\proves\theta\}$}
is recursive.
Show that any complete recursively axiomatizable theory is decidable.
Show that the theory of $(\omega,Sc,0)$ is decidable. Show that
the theory of dense linear orderings with no end points is decidable.
Show that the theory of $({\Bbb R},E)$ is decidable, where
$E=\{(x,y):(x-y)\in {\Bbb Q}\}$.  Show that the $Th({\Bbb C},+,\cdot)$
is decidable.

\prob Show that any undecidable recursively axiomatizable theory is
incomplete.

\prob Let $T$ be a theory in a finite language without any infinite
models.  Show that $T$ is decidable.

\prob Suppose that $T$ is a recursively axiomatizable theory.
Show that $T$ is decidable iff the set of sentences $\theta$
such that $T\cup\{\theta\}$ has a model is recursively enumerable.

\prob Show that the set of validities in the language of pure equality
is decidable.

\com{Show there exists a recursive list of all complete theories.}

\prob Show that the theory of one equivalence relation is decidable.

\prob  Let $T$ be a consistent recursively axiomatizable theory.
Call a formula $\theta(x)$ in the language of $T$ strongly finite
iff in every model $\moda$ of $T$, only a finite number of elements
satisfy $\theta(x)$.  Prove that the set of strongly finite
formulas is recursively enumerable.

\prob Suppose $T$ is a decidable theory and $\theta$ is a sentence
in the language of $T$.  Show that $T\cup\{\theta\}$ is decidable.
Hence finite extensions of decidable theories are decidable.

\prob Suppose $T$ is a consistent decidable theory.
Show there exists a complete consistent decidable $T^\prime\supseteq T$ in
the same language as $T$.

\prob Suppose $T$ is a consistent decidable theory in a language $L$.
Suppose that $L^{\prime}=L\cup\{c\}$ is the bigger language
with one new constant symbol.  Prove that the $L^{\prime}$ theory
axiomatized by $T$ is decidable.

\prob Suppose $T$ is a consistent decidable theory.
Show there exists a consistent Henkin $T^\prime\supseteq T$
which is decidable.

\prob Suppose $T$ is a consistent decidable theory.
Show that $T$ has a \dex{recursive model} (i.e.,  a model whose
universe is a recursive subset of $\omega$ and such that all
relations and operations are recursive).
\sector{Undecidable Theories}

The \dex{incompleteness theorem} was proved by
Kurt G\"{o}del in 1930.  It was improved by Barkley Rosser in 1936.
A vague statement of it is that any
theory which can be recursively  axiomatized  and is strong enough to
prove elementary facts about multiplication and addition on the
natural numbers must be incomplete.  Here we will develop
a version of this theorem based on a small fragment of the theory of
finite sets.

\dex{HF} is the set of all {hereditarily finite sets}.
A set is \dex{hereditarily finite} iff
the set is finite,
all elements of it are finite, all elements of elements of it are finite,
etc, etc.  The formal definition is that $HF=\cup_{n<\omega}V_n$
where $V_0=\emptyset$ and $V_{n+1}$ is the set of all subsets of $V_n$.

\prob
Show that if $x,y\in HF$ then $\langle x,y\rangle\in HF$.
Show that if $X,Y\in HF$ then $X \times Y \in HF$.
Let ${Y^X}$ be the set of all functions with domain $X$ and
range $Y$.  Show that if $X,Y\in HF$ then $Y^X \in HF$.

\prob
Prove that $HF$ is transitive.

\bigskip

The ${\Delta_0}$ formulas \sdex{$?Delta_0$} are the smallest
set of formulas containing
$x \in y$, $x=y$ (i.e., all atomic formulas),
closed under $\disj$ and $\neg$ and
bounded quantifiers $\forall x\in y$ and $\exists x\in y$.
By this we mean that if $\theta$ is a $\Delta_0$ formula and
$x$ and $y$ are any two variables, then
$\exists x\in y\;\theta$ and $\forall x\in y\;\theta$ are $\Delta_0$ formulas.
The formal definition of $\exists x\in y\;\theta$ is
$\exists x (x\in y\conj\theta)$ and
$\forall x\in y\;\theta$ written formally is
$\forall x(x\notin y\disj\theta)$.

A formula is ${\Sigma_1}$   iff it has the form:
$$\exists x_1 \exists x_2\ldots \exists x_n \phi$$
where $\phi$ is a $\Delta_0$ formula. \sdex{$?Sigma_1$}
A subset  $X$ of $HF$ is $\Delta_0$
iff there exists a formula $\theta(x)$ (which may have constant symbols for
elements of $HF$) and $X=\{a\in HF: \HFin\models\theta(a)\}$.  Similarly for
$\Sigma_1$ subsets.  A subset of $HF$ is ${\Delta_1}$
\sdex{$?Delta_1$}
iff it and
its complement are $\Sigma_1$.

\bigskip

\prob
Show that $\omega$ is a definable subset of $\HFin$.
Show that $\omega$ is a $\Delta_1$ subset of $\HFin$.

\prob $A^{<\omega}$ is the set of all functions with domain some
$n\in\omega$ and range contained in $A$. Show for
any $A\in HF$ that $A^{<\omega}$ is $\Delta_1$.

\prob Suppose $\partial m:S\times A\to S\times A\times D$ is a
partial function
in $HF$ and $B\subseteq A$ and $f:B^{<\omega}\to B^{<\omega}$
is the associated
partial recursive function.  Show there exists a
$\Sigma_1$ formula $\theta(x,y)$ such that for every
$u,v\in B^{<\omega}$
$$f(u)=v\mbox{ iff } \HFin\models\theta(u,v)$$

\prob Show that a set $X\subseteq B^{<\omega}$ is recursively enumerable
iff it is a $\Sigma_1$ subset of $HF$.

\prob Show that for any $B\in HF$ every recursive subset of $B^{<\omega}$
is a
$\Delta_1$ subset of $HF$.  Show that every $\Delta_1$ subset of
$B^{<\omega}$ is recursive.

\prob  Assume that every recursively enumerable subset of $\omega$ is
a definable subset of $HF$.
Show that the complete theory $Th(HF,\in)$ is undecidable.
Show that every recursively axiomatizable $T\subseteq Th(HF,\in)$
is incomplete.

\prob
Consider the following theory of finite sets \dex{FIN}.  The language
of FIN has only one binary relation $E$ which we like to write
using infix notation. There are four axioms:

 Empty Set. $ \exists x\, \forall y \neg ( y E x )$

 Extensionality.
 $ \forall x \forall y ( x = y \iff \forall z ( z E x \iff z E y)) $

 Pairing. $ \forall x \forall y \exists z \forall u (
u E z \iff (u = x \disj u = y))$

 Union. $ \forall x\, \exists y\, ( \forall z ( z E y \iff
(\exists u\, u E x \conj z E u ))$

\noindent Prove that $\HFin$ is a model of FIN.

\prob The theory FIN has no constant symbols,  however the
empty set is unique and can be
defined by $\theta(x)=\forall y  \neg (yEx)$.   Show there
exists an effective mapping from $HF$ into $\Delta_0$ formulas
(without parameters and having one free variable $x$),
say $a\to \theta_a(x)$, such that
$$\HFin\models\theta_a(a)$$ and
$$\mbox{FIN }\proves\exists ! x\;\theta_a(x).$$

\prob
If $M$ and $N$ are models of FIN we say that ${M\subseteq_e N}$
\sdex{$M?subseteq_e N$}
 iff
$M$ is a substructure of $N$ and $N$ is an \dex{end extension} of
$M$, this means that for any $a\in M$ if $N\models b E a$ then
$b\in M$.
 Show that for every model $N$ of FIN there exists $M\subseteq_e N$
such that $M$ is isomorphic to $\HFin$.

\prob
Let $\theta(u)$ be any $\Sigma_1$ formula.  Show that for any
$u\in HF$:
$$\HFin\models\theta(u)\mbox{ iff FIN}\proves\theta(u)$$
Since FIN has no constant symbols we really mean
$\exists x(\theta_u(x)\conj\theta(x))$ where we have written
$\theta(u)$.

\prob Show that FIN is not a decidable theory, i.e.,
the set
$$\{\theta: \mbox{FIN}\proves\theta\}$$
is not recursive.

\prob (Church) By considering validities of the
form $\conj\mbox{FIN}\implies\theta$,
show that the set of validities in the language of
one binary relation is not recursive.

\prob Let $T$ be any consistent recursively axiomatizable theory in one
binary relation that is consistent with FIN (eg ZFC).  Show that
$T$ is incomplete.

\com{This asks only to show incompleteness without explicitly giving
the sentence which is not decided by $T$.
But note that if $x\in K$ iff $T\proves\theta(x)$,
where $K=\{e:e\in W_e\}$, then let $W_e=\{n:T\proves\neg\theta(n)\}$.
Then we get that $\theta(e)$ is not decided by $T$.}

\prob (Mostowski, Robinson, Tarski)
Call a theory $T$ \dex{strongly undecidable} iff $T$ is finitely
axiomatizable and every theory
$T^\prime$ in the language of $T$ that is consistent with $T$
is undecidable.  Show that FIN is strongly undecidable.

\prob We say that a model $\moda$ is \dex{interpretable}
in $\modb$ iff the universe of $\moda$ and the relations and operations
of $\moda$ are definable in $\modb$. They may have completely different
languages, say $L_a$ and $L_b$.  Assume these languages are finite.
Show there exists a recursive map $r$ from $L_a$ sentences to
$L_b$ such that
$$\moda\models\theta \rmiff \modb\models r(\theta)$$
for any $L_a$ sentence $\theta$.
This technique is called relativization of quantifiers.

\prob  Call a structure \dex{strongly undecidable} iff it models some
strongly undecidable theory.
Suppose $\modb$ is a structure in which a strongly undecidable
structure $\moda$ is definable.
Show that $\modb$ is strongly undecidable.

\prob Show the theory of  graphs is undecidable, i.e., the set
of sentences true in every graph is not recursive.
Hint:
Let $(A,R)$ be any binary relation and find a graph
$(V,E)$ and a pair of formulas $U(x)$ and $Q(x,y)$
 in the language of graph theory such that
$A=\{v\in V:(V,E)\models U(v)\}$ and $R=\{(u,v):(V,E)\models Q(u,v)\}$.

\prob Show that the theory of partially ordered sets is undecidable.

\com{
Let $(V,E)$ be any graph and consider the partial order
$(V\cup E,\leq)$ where $a\leq b$ iff $a=b$ or for some $x,y\in V$
$a=x$ and $b=(x,y)\in E$.
}

\prob Show that the theory of two unary operation symbols is
undecidable.

\com{
Given a graph $(V,E)$ consider the structure
$(V\cup E,f,g)$ where $f((a,b))=a$ and $g((a,b))=b$ and $f(a)=g(a)=a$.
}

\prob Show that the structure $(\omega,+,\cdot,0,1,\leq, x^y)$ is strongly
undecidable.
Hint:
Use only positive integers.
Define $xEy$ iff there exists a prime $p$ such that $p^x$ divides
$y$ but $p^{x+1}$ does not.  Say that $z$ is transitive iff
for any $x,y$ if $xEy$ and $yEz$ then $xEz$.
Define the $x\eq y$ iff for all $z$
$z E x$ iff $z E y$.  Define $x$ to be minimal iff for all
$y\eq x$ we have $x\leq y$. Let $H$ be the set of all $x$ such there
exists a transitive $z$ with $xEz$ and every $y$ such that $yEz$
is minimal.
Show that there exists an isomorphism from $(H,E)$ to $\HFin$.

\com{
Let $f:HF\to H$ be defined by induction on the
rank. $f(\emptyset)=1$. suppose $f(a)$ defined for every
$a\in b$, then consider the smallest number
$z$ such that for every $a\in b$ there exists a prime $p$
such that  $p^{f(a)}$ divides $z$.  Then $z\in H$ and
we let $f(b)=z$.
}

\bigskip

The next goal is show that the structure $(\omega,+,\cdot)$ is strongly
undecidable.   To do this we need a little number theory.

\prob Prove the \dex{Chinese Remainder} Theorem.
Given any finite sequence
of pairwise relatively prime integers $\langle p_0,\ldots,p_{n-1}\rangle$
and given any finite sequence of integers $\langle x_0,\ldots,x_{n-1}\rangle$
there exists an integer $x$ such that $x_i=x$ mod $p_i$ for all $i<n$.
Two integers are \dex{relatively prime}
iff their greatest common divisor is one.

\com{
Let $q_i$ be the product of all $p_j$ for $j\not=i$.  Then $q_i$ and
$p_i$ are relatively prime, so there exists $u_i$ with
$u_iq_i=1$ mod $p_i$. Now let $x=\Sigma x_iu_iq_i$.
}

\prob Prove that $p_i=1+(i+1)(n!)$ for $i<n$ are
pairwise relatively prime.

\com{If $x$ divides $1+in!$ and $1+jn!$ then it divides
$(i-j)n!$, so the smallest prime $x$ is less than $n$, contradicting
$x$ divides $1+in!$}

\prob Show there exists a formula $\theta(i,u,x,y)$
in the language of $(\omega,+,\cdot)$ such
that for every $\langle x_0,\ldots,x_{n-1}\rangle\in\omega^{<\omega}$
there exists $x,y\in\omega$ such that for each $i<n$
$x_i$ is the unique $u$ such $(\omega,+,\cdot)\models \theta(i,u,x,y)$.

\prob Show that the structure $(\omega,+,\cdot)$ is strongly
undecidable.

\com{
The only difficulty is to find a definition of $z=x^y$.
Change the recursive definition ($x^0=1$ and $x^{y+1}=x^yx$)
into an explicit definition by using $\theta$.
}

\prob Show that the structure $({\Bbb Z},+,\cdot)$ is strongly
undecidable.
Hint:Lagrange proved that every positive integer is the sum of four
squares.

\prob Prove that the theory of rings is undecidable.

\prob Show that the structure $({\Bbb Z},+,x$ divides $y)$ is
strongly undecidable.

\com{
You need to define multiplication.  To define $x^2$ note that
it is the least integer $y$ such $x$ and $x+1$ both divide
$y+x$.   To define $xy$ note that
$(x+y)^2=x^2+xy +xy+y^2$.
}

\prob (Tarski) Show that the theory of groups is undecidable.
Hint: Let $G$ be the group of all bijections of ${\Bbb Z}$ into itself
where the group operation is composition. Let $s\in G$ be defined
by $s(x)=x+1$.  Embed ${\Bbb Z}$ into $G$ by mapping
$i$ to $s^i$.  Show that for any $t\in G$ that $t$ commutes
with $s$ iff $t=s^i$ for some integer $i$.  Show
that $i$ divides $j$ iff every element of $G$ that commutes with
$s^i$ commutes with $s^j$.
Show therefore that $({\Bbb Z},+,i$ divides $j)$ can
be defined in $G$.

\com{ If $f$ commutes with $s$ then $f(i+1)=f(i)+1$ so
easily $f(0+x)=f(0)+x=s^{f(0)}(x)$. To say that $f$ commutes
with $s^i$ means that $f(x+i)=f(x)+i$ for every $x$.  So
for example, let $f(x)=x+i$ for $x$ a multiple of $i$ and
$f(x)=x$ otherwise, then if $i$ does not divide $j$ then
$j=f(0+j)\not=f(0)+j=i+j$, so $f$ does not commute with
$s^j$.
Actually we are just getting that $(G,s)$ is strongly undecidable,
but probably we get $G$ strongly undecidable, by existentially quantifying
out $s$ from the finite theory?}

\prob Julia Robinson showed that ${\Bbb Z}$ is
explicitly definable in the structure $({\Bbb Q},+,\cdot)$.
Prove that the theory of fields is undecidable.

\prob Show that ${\Bbb Z}$ is not definable in $({\Bbb C},+,\cdot)$.

\prob (Julia Robinson)
Show that plus is definable in $(\omega, \cdot,Sc)$.
Hence the structure $(\omega, \cdot,Sc)$ is strongly undecidable.
Hint: Let $a+b=c$ and multiply out $(ac+1)(bc+1)$.

\prob Prove or give a counterexample:  Let $S$ be a decidable
set of sentences in propositional logic using the propositional
letters $\{P_n:n<\omega\}$. Then
$$\{\theta: S\proves\theta\}$$
is decidable.

\prob Let $T$ be a recursively axiomatizable theory with no decidable
consistent complete extensions.   Show that there are $\cc$ many distinct
complete extensions of $T$.

\prob Let $T$ be a consistent recursively axiomatizable extension
of FIN.
Let $A\subseteq\omega$ be a nonrecursive set of integers.
Show there exists a model $\moda$ of $T$ that is an end extension
of $\HFin$ and such that
for no formula $\theta(x)$ is
$$A=\{n<\omega:\moda\models\theta(n)\}.$$

\com{For any $T^\prime$  a consistent recursively axiomatizable extension
of FIN and formula $\theta(x)$ the sets
$\{n:T^\prime\proves\theta(n)\}$ and $\{n:T^\prime\proves\neg\theta(n)\}$ are
disjoint recursively enumerable sets.}

\prob Prove that the following set is recursively enumerable
but not recursive:
$$\{\theta :\theta \mbox{ is an } L\mbox{ sentence with a finite model}\}$$

\com{Think of a computation as a finite model of some theory.}

\prob
Another way to see that the structure $(\omega,+,\cdot,0,1,\leq)$
is strongly undecidable is to use Robinson's theory $Q$.
The axioms of ${Q}$ are \sdex{$Q$}
\begin{enumerate}
 \item $0\not=1$
 \item $\forall x (x+0=x)$
 \item $\forall x\forall y(x+(y+1)=(x+y)+1$
 \item $\forall x(x\cdot 0=0)$
 \item $\forall x\forall y(x\cdot(y+1)=x\cdot y +1)$
 \item $\forall x\forall y(x\leq y\iff\exists z(x+z=y))$
\end{enumerate}
Then the $\Delta_0$ are the smallest family of formulas containing
the atomic formulas, closed
under $\disj,\neg$ and bounded quantification, i.e., $\exists x\leq y$
and $\forall x\leq y$.   Show that
a set $A\subseteq \omega$ is recursively enumerable iff there
exists a $\Sigma_1$ formula $\theta(x)$ such that for $n<\omega$
$$n\in A\mbox{ iff } Q\proves \theta(n)$$

\prob The theory of \dex{Peano arithmetic}, \dex{PA}, is
axiomatized by $Q$ plus the induction axioms.  For each formula
$\theta(x,\vec{y})$ the following is an axiom of PA:
$$\forall\vec{y}[(\theta(0,\vec{y})\conj\forall x (\theta(x,\vec{y})\implies
\theta(x+1,\vec{y})))\implies \forall x \theta(x,\vec{y})]$$
The analogue of PA in set theory is ZFCFIN.  \dex{ZFCFIN} is ZFC minus
the axiom of infinity plus the negation of the axiom of infinity.
PA and ZFCFIN have exactly the same strength.  There is a correspondence
between models and theorems of these two theories.
Formalize these statements and prove them.

\sector{The second incompleteness theorem}

G\"{o}del's \dex{second incompleteness} theorem may be paraphrased by
saying that a consistent recursively axiomatizable theory
that is sufficiently
strong cannot prove its own consistency unless it is inconsistent.
To prove it we must give a more explicit proof
of the first incompleteness theorem.

Thruout this section we will be assuming
$T$ is a consistent recursively axiomatizable theory in a
language that contains a symbol $n$ for each $n\in\omega$.
(More generally it would be enough to suppose these constants were
nicely definable in $T$, as for example they are in set theory.)

We say that the partial recursive functions are \dex{representable} in $T$ iff
for every partial recursive function $f:\omega\mapsto\omega$
there exists a formula $\theta(x,y)$ such
that for every $n,m\in\omega$
$$f(n)=m\mbox{ iff }T\proves\theta(n,m)$$
and furthermore
$$T\proves\forall x\exists^{\leq 1}y\;\theta(x,y)$$
the quantifier $\exists^{\leq 1}y$ meaning there exists at most one $y$.
Note that we do not demand that $f$ is total.

We say that the recursively enumerable sets are representable in $T$ iff
for every recursively enumerable set $A\subseteq \omega$ there
is a formula of $\theta(x)$ such that for
every $n\in\omega$,
$$n\in A \mbox{ iff } T\proves\theta(n)$$
We have already seen that a consistent recursively axiomatizable theory
in which every recursively enumerable set is representable must be
incomplete, since there exist recursively enumerable sets that are
not recursive.

\prob Show that if every partial recursive function is representable
in $T$ then every recursively enumerable set is representable in $T$.
\com{I don't know if the converse is true.}

\bigskip
Now since $T$ is recursively axiomatizable the formulas of the language
of $T$ are a recursive set.
We use $\theta\mapsto \ul\theta \ur$ to mean
a recursive map from the sentences of the language
of $T$ to $\omega$. \sdex{$?ul?theta ?ur$}
One could get tricky and do this in some unnatural
way, but we assume here that is done simply, as follows.
First number all the symbols of $T$ and the logical symbols.
We assume that the type and
arity of each symbol is given by recursive functions and
also our special symbols `$n$' correspond to a recursive set.
The map that takes $n$ to the number that codes the symbol
$n$ should be recursive.
Then each formula is a finite string of symbols and so corresponds
to a sequence of numbers
$\langle a_1,a_2,\ldots, a_n\rangle$.  Map this to
$2^{a_1}3^{a_2}\ldots p_n^{a_n}$ where $p_n$ is the $n^{th}$ prime
number.   Such a map (from sentences to numbers)
is called a G\"{o}del numbering.

\bigskip

{\bf Fixed point lemma:}
Suppose $T$ is a consistent recursively axiomatizable theory
in which every partial recursive function is representable.
Let $\psi(x)$ be any formula with one free variable.
Then there exists a sentence $\theta$ such that
$$T\proves\mbox{`} \theta \iff\psi(\ul\theta\ur)\mbox{'}$$

\proof
The sentence $\theta$ says in effect
``I have the property defined by $\psi$''.
Consider the recursive map
$\rho(x)\mapsto \rho(\ul \rho \ur)$.  That is, given a formula
$\rho(x)$ with one free variable map it to the formula gotten
by substituting the G\"{o}del number of $\rho$ into the free variable
of $\rho$.  Since this function is recursive
there exist a formula $\chi(x,y)$ such that
for any $\tau$ and $\rho(x)$ we have that
$$\tau =\rho(\ul \rho \ur) \mbox{ iff }
T\proves \chi(\ul\tau\ur,\ul\rho\ur)$$
and
$$T\proves\forall x\exists^{\leq 1}y\;\chi(x,y)$$
Now define $\sigma(x)=\exists y (\chi(x,y)\conj \psi(y))$ and
let $\theta=\sigma(\ul\sigma\ur)$.   We claim that
$T\proves \psi(\ul\theta\ur)\iff\theta$.
To see this let $\moda$ be any model of $T$.  If $\moda\models\theta$ then
$\moda\models\exists y (\chi(\ul\sigma\ur,y)\conj \psi(y))$, by
the definition of $\theta$.  But $T\proves \chi(\ul\sigma\ur,\ul\theta\ur)$,
hence $\moda\models\psi(\ul\theta\ur)$.  Alternatively
suppose $\moda\models\psi(\ul\theta\ur)$ then
$\moda\models\chi(\ul\sigma\ur,\ul\theta\ur)\conj\psi(\ul\theta\ur)$
Hence $\moda\models\exists y \;\chi(\ul\sigma\ur,y)\conj\psi(\ul\theta y\ur)$,
and so $\moda\models\theta$.
\qed

\prob (Tarski, truth is not definable.)
Suppose $T$ is  a consistent recursively axiomatizable theory
in which every partial recursive function is representable and
$\moda$ is a model of $T$.
Let Truth$=\{\ul\theta\ur:\moda\models\theta\}$.
Show that Truth is not definable in $\moda$.

\com{Let $\psi(x)$ be such a definition and apply the fixed point
theorem to $\neg\psi(x)$.}

\prob  (The First Incompleteness Theorem.)
If $T$ is  a consistent recursively axiomatizable theory
in which every partial recursive function is representable, then
the theorems of $T$ are a recursively enumerable set.
Hence there exists
a formula $\sdex{?prf(x)}$ in the language of $T$ such that
for any $\theta$ a sentence in the language of $T$
$$T\proves \prf({\ul \theta \ur})\mbox{ iff }
T\proves \theta $$
Show there exists a sentence $\theta$
that asserts its own unprovability, that is
$$T\proves \mbox{`}\theta\iff \neg \prf(\ul \theta \ur)\mbox{'}$$
Show that this $\theta$ is neither provable nor refutable in $T$.

\com{ Consider $\psi(x)=\neg \prf(x)$.
If $T\proves\theta$, then
$T\proves \prf(\ul \theta \ur)$ by representability,  but
by the theorem, $T\proves\neg \prf(\ul \theta \ur)$
and then $T$ is inconsistent.  If $T\proves\neg\theta$, then
$T\proves \prf(\ul \theta \ur)$ by the theorem, and
hence by representability $T\proves \theta$, so $T$ is inconsistent.}

\bigskip

{\bf The Second Incompleteness Theorem:}
From now on
fix a consistent recursively axiomatizable theory $T$
in which every partial recursive function is representable.
The statement $\sdex{con(T)}$ which stands for $T$ is consistent
can now be thought of as a sentence in the language of $T$:
$$con(T)=\neg \prf(\ul 0=1 \ur)$$
We wish to show that $T$ does not prove $con(T)$.  However
there are many different $\prf$ formulas and for some of
these it may be possible to prove $con(T)$!

\prob Let $\prf^*(x)$ be the formula
$(\prf(x)\conj x\not=\ul 0=1 \ur)$. Show
that $T\proves con^*(T)$ where $con^*(T)=\neg \prf^*(\ul 0=1\ur)$,
and for every sentence $\theta$ we have
$T\proves\theta$ iff $T\proves \prf^*(\ul \theta \ur)$.

We say that the $\prf(x)$ is a \dex{reasonable proof predicate} for a theory
$T$ iff it satisfies:
\begin{enumerate}
  \item $T\proves \prf({\ul \theta \ur})\mbox{ iff }
         T\proves \theta$.
  \item  $T\proves$`$(\prf(\ul\theta\ur)\conj \prf(\ul(\theta\implies\rho)\ur)
         \implies \prf(\ul\rho\ur)$'.
  \item  $T\proves$`$\prf(\ul\theta\ur)\implies \prf(\ul
          \prf(\ul \theta \ur)\ur)$'.
\end{enumerate}
The second condition is kind of a Modus ponens for proofs while the
last condition basically says that if you can prove $\theta$, then
you can prove that you can prove $\theta$.

\prob Show that \zfc $\;\;$ has a reasonable proof predicate.

\prob Suppose that $T$ is a theory with a reasonable proof
predicate $\prf(x)$ and $\theta$ is a sentence which asserts its
own unprovability:
$$T\proves\theta\iff \neg \prf(\ul\theta\ur)$$
Show that $T\proves con(T)\implies\theta$.  Conclude that
$T\not\proves con(T)$.

\com{ Show:
a. $T\proves \prf(\ul\neg\theta\ur)\iff \prf(\ul \prf(\ul\theta\ur)\ur)$.
b. $T\proves (\prf(\ul\theta\ur)\conj \prf(\ul\neg\theta\ur))
\implies \prf(\ul 0=1\ur)$.
c. $T\proves \prf(\ul\theta\ur)\implies \prf(\ul
          \prf(\ul \theta \ur)\ur)$
Combining a. and c. we get
$T\proves \prf(\ul\theta\ur)\implies \prf(\ul\neg\theta\ur))$.
Using b. we get $T\proves \prf(\ul\theta\ur)\implies \prf(\ul 0=1\ur)$.
Taking the contrapositive
So $T\proves con(T)\implies\neg \prf(\ul\theta\ur)$.
But $T\proves \neg \prf(\ul\theta\ur)\implies\theta$.  Hence
$T\proves con(T)\implies\theta$. }

\prob (L\"{o}b)  Let $\theta$ be a sentence of \zfc such that
$$\mbox{\zfc}\proves \mbox{`}\theta\iff \prf(\ul
\theta \ur)\mbox{'}$$
Show that \zfc$\proves\theta$.

\com{
Note that \zfc$\proves con(\zfc+\neg\theta)\iff\neg
\prf(\ul \theta \ur)$.
Since $$\zfc\proves\theta\iff \prf(\ul\theta\ur),$$ by putting these together
 we get \zfc+$\neg\theta\proves con(\zfc+\neg\theta)$.
But this is a theory that proves its own consistency and hence
by the second incompleteness theorem it must be inconsistent,
which means \zfc$\proves\theta$.
}

\prob Suppose \zfc is consistent.  Show that there exists a recursively
axiomatizable consistent $T$ extending \zfc such that
$T\proves \neg con(T)$.

\com{T is zfc+$\neg$con(zfc). Need to see that con(T) implies con(zfc)
This is clear for the reasonable proof predicate that we constructed.}

\prob Find a consistent theory $T\supseteq \zfc$ and a first order sentence
(in any language) such that
\par\centerline{$T\proves\theta$ has a finite model}\par\noindent
but $\theta$ does not have a finite model.

\com{same $T$ works.  Models of $\theta$ are Turing machine computations
that enumerate a proof of $\neg con(\zfc)$.  This can be done also
with just representability of $K$: if $R(x,y)$ is recursive
and $K$ is the x projection of $R$,
for some $n\notin K$ \zfc$+\exists y\in\omega\; R(n,y)$ is consistent, (but
there can't be a finite $y$).
We can take $\theta$ to give
`topped' models of number theory.}

\prob Suppose $T$ is a recursively axiomatizable theory
contained in $Th(\omega,+,\cdot)$. Show that there exists a
model $\moda$ of $T$ and a formula $\theta(x)$ such that
the least $n\in A$ such that $\moda\models\theta(n)$ is infinite.

\com{This can be done by constructing a reasonable proof predicate
and letting $n$ be a proof of the inconsistency of FIN.
It can also be done by simply finding $R$ recursive
that $T\proves\neg R(n)$, for all finite $n$ but
$T$ is consistent with $\exists x R(x)$.}
\sector{Further Reading}

\prob Mathematical Logic, J.R.Shoenfield, Addison-Wesley, 1967.

\prob Handbook of Mathematical Logic, edited by Jon Barwise,
   Studies in Logic and the Foundations of Mathematics, vol.90,
   North Holland, 1977.

\prob Model Theory, C.C.Chang and H.J.Keisler,
   Studies in Logic and the Foundations of Mathematics, vol.73,
   North Holland, 1973.

\prob Some aspects of the Theory of Models, R.L.Vaught,
   in Papers in the Foundations of Mathematics, Number 13 of
   the Herbert Ellsworth Slaught Memorial Papers, Mathematics
   Association of America, 1973, 3-37.

\prob Set Theory, K.Kunen,
   Studies in Logic and the Foundations of Mathematics, vol.102,
   North Holland, 1980.

\prob Set Theory, T.Jech, Academic Press, 1978.

\prob Recursively Enumerable Sets and Degrees, R.I.Soare,
   Perspectives in Mathematical Logic, Springer-Verlag, 1987.

\prob Self-Reference and Modal Logic, C.Smory\'{n}ski,
    Springer-Verlag, 1985.

\prob Undecidable Theories, A.Tarski, A.Mostowski, and R.M.Robinson,
   North-Holland, 1953.

\prob  Journal of Symbolic Logic, Annals of Pure and Applied Logic,
   Fundamenta Mathematicae.

\prob Electronic preprints in Logic:
email listserv@math.ufl.edu with  ``HELP settheory'' in the
body of the message.  Or you can
anonymously ftp to
$$\mbox{ftp.math.ufl.edu}$$
and look in directory /pub/logic.

\newpage


%
%

 \twocolumn\pagestyle{myheadings}           
 \markboth{Index}{Index}\markright{Index}   
    \par\noindent $(A,R)\isom (B,S)$ \dotfill 3.4-p.15
    \par\noindent $(\moda,a)_{a\in A}$ \dotfill 7.1-p.29
    \par\noindent $A \cup B$ \dotfill 2.1-p.10
    \par\noindent $cf(\beta)$ \dotfill 6.6-p.22
    \par\noindent $D(\moda)$ \dotfill 8.1-p.44
    \par\noindent $D_0(\moda)$ \dotfill 8.23-p.46
    \par\noindent $e\res {\pp}$ \dotfill 1.3-p.4
    \par\noindent $f\res C$ \dotfill 2.6-p.11
    \par\noindent $f^{\prime \prime }C$ \dotfill 2.6-p.11
    \par\noindent $K$ \dotfill 9.23-p.55
    \par\noindent $L$ theory \dotfill 7.2-p.36
    \par\noindent $L_{\moda}=L\cup\{c_a: a\in A\}$ \dotfill 7.1-p.29
    \par\noindent $M\subseteq_e N$ \dotfill 11.12-p.63
    \par\noindent $Q$ \dotfill 11.39-p.66
    \par\noindent $Th(\moda)$ \dotfill 7.5-p.37
    \par\noindent $V$ \dotfill 2.1-p.11
    \par\noindent $Val$ \dotfill 10.11-p.59
    \par\noindent $V_{\alpha }$ \dotfill 5.14-p.18
    \par\noindent $x + 1 = x \cup \{x\}$ \dotfill 2.1-p.10
    \par\noindent $X \times Y$ \dotfill 2.4-p.11
    \par\noindent $X\cap Y$ \dotfill 2.2-p.11
    \par\noindent $X\setminus Y$ \dotfill 2.2-p.11
    \par\noindent $y = P(x)$ \dotfill 2.1-p.10
    \par\noindent $y = \cup x $ \dotfill 2.1-p.10
    \par\noindent $y = \{ x \in z : \theta(x) \}$ \dotfill 2.1-p.10
    \par\noindent $Y^X$ \dotfill 2.5-p.11
    \par\noindent $z = \{ x,y \}$ \dotfill 2.1-p.10
    \par\noindent $z \subseteq x $ \dotfill 2.1-p.10
    \par\noindent $[X]^{<\omega}$ \dotfill 2.23-p.13
    \par\noindent $[X]^{\omega}$ \dotfill 2.31-p.13
    \par\noindent $[\kappa ]^{<\omega }$ \dotfill 6.13-p.23
    \par\noindent $[\omega]^{<\omega}$ \dotfill 2.23-p.13
    \par\noindent $[\omega]^{\omega}$ \dotfill 2.29-p.13
    \par\noindent $\aleph _{\alpha }$ \dotfill 6.1-p.22
    \par\noindent $\bigcap X$ \dotfill 2.2-p.11
    \par\noindent $\cc$ \dotfill 2.18-p.12
    \par\noindent $\Delta_0$ \dotfill 11.3-p.62
    \par\noindent $\Delta_1$ \dotfill 11.3-p.62
    \par\noindent $\emptyset$ \dotfill 2.1-p.10
    \par\noindent $\kappa +\gamma $ \dotfill 6.3-p.22
    \par\noindent $\kappa \gamma $ \dotfill 6.3-p.22
    \par\noindent $\kappa ^+$ \dotfill 6.11-p.22
    \par\noindent $\kappa ^{<\gamma }$ \dotfill 6.12-p.23
    \par\noindent $\kappa ^{\gamma }$ \dotfill 6.12-p.23
    \par\noindent $\langle x,y\rangle$ \dotfill 2.3-p.11
    \par\noindent $\moda\elemsub^j\modb$ \dotfill 8.1-p.43
    \par\noindent $\moda\equiv\modb$ \dotfill 7.6-p.37
    \par\noindent $\moda\models\theta$ \dotfill 7.1-p.28
    \par\noindent $\moda\models\theta(a_1,\ldots,a_n)$ \dotfill 8.1-p.40
    \par\noindent $\moda\res{L_1}$ \dotfill 7.1-p.29
    \par\noindent $\modn_{Sc}$ \dotfill 8.8-p.45
    \par\noindent $\omega $ \dotfill 5.1-p.17
    \par\noindent $\omega _{\alpha }$ \dotfill 6.1-p.22
    \par\noindent $\omega$ \dotfill 2.1-p.10
    \par\noindent $\ord$ \dotfill 5.9-p.18
    \par\noindent $\Sigma\proves\theta$ \dotfill 10.1-p.56
    \par\noindent $\Sigma_1$ \dotfill 11.3-p.62
    \par\noindent $\ul\theta \ur$ \dotfill 12.2-p.68
    \par\noindent ${\cal M}(T)$ \dotfill 7.10-p.37
    \par\noindent ${\moda} \elemsub {\modb}$ \dotfill 8.1-p.40
    \par\noindent ${\moda} \isom {\modb}$ \dotfill 8.1-p.42
    \par\noindent ${\moda} \substruc {\modb}$ \dotfill 8.1-p.40
    \par\noindent ${\qq}$ \dotfill 2.25-p.13
    \par\noindent ${\rr}$ \dotfill 2.18-p.12
    \par\noindent ${\zz}$ \dotfill 2.19-p.13
    \par\noindent $| X | < |Y |$ \dotfill 2.11-p.12
    \par\noindent $| X | = | Y |$ \dotfill 2.11-p.12
    \par\noindent $| X | \leq | Y |$ \dotfill 2.11-p.12
    \par\noindent $|{\moda}|$ \dotfill 8.1-p.42
    \par\noindent adequate \dotfill 1.1-p.4
    \par\noindent algebraic \dotfill 2.26-p.13
    \par\noindent almost disjoint family \dotfill 2.36-p.13
    \par\noindent arity \dotfill 7.1-p.25
    \par\noindent atomic diagram \dotfill 8.23-p.46
    \par\noindent atomic formulas \dotfill 7.1-p.25
    \par\noindent atomic sentences \dotfill 1.1-p.4
    \par\noindent Axiom of Choice \dotfill 4.1-p.16
    \par\noindent axiomatization \dotfill 7.2-p.36
    \par\noindent axiomatizes \dotfill 7.2-p.36
    \par\noindent binary relation \dotfill 1.24-p.7
    \par\noindent bipartite graph \dotfill 2.10-p.12
    \par\noindent bound \dotfill 7.1-p.26
    \par\noindent canonical structure \dotfill 7.1-p.33
    \par\noindent Cantor normal form \dotfill 5.21-p.20
    \par\noindent cardinal \dotfill 6.1-p.22
    \par\noindent cartesian product \dotfill 2.4-p.11
    \par\noindent categorical \dotfill 8.7-p.45
    \par\noindent chain \dotfill 1.21-p.6
    \par\noindent Chinese Remainder \dotfill 11.24-p.65
    \par\noindent chromatic number \dotfill 1.30-p.8
    \par\noindent Church's Thesis \dotfill 9.22-p.52
    \par\noindent class \dotfill 2.1-p.10
    \par\noindent cofinality \dotfill 6.6-p.22
    \par\noindent cofinite \dotfill 8.10-p.45
    \par\noindent complete \dotfill 1.15-p.6
    \par\noindent complete \dotfill 7.1-p.31
    \par\noindent components \dotfill 2.9-p.11
    \par\noindent con(T) \dotfill 12.4-p.70
    \par\noindent connected \dotfill 2.9-p.11
    \par\noindent consistent \dotfill 10.1-p.57
    \par\noindent consistent \dotfill 7.19-p.38
    \par\noindent contradiction \dotfill 1.1-p.4
    \par\noindent countable \dotfill 2.18-p.12
    \par\noindent decidable \dotfill 10.14-p.60
    \par\noindent definable \dotfill 8.6-p.45
    \par\noindent dimension \dotfill 1.29-p.8
    \par\noindent disjunctive normal form \dotfill 1.15-p.5
    \par\noindent elementarily equivalent \dotfill 7.6-p.37
    \par\noindent Elementary Chain Lemma \dotfill 8.28-p.47
    \par\noindent elementary diagram \dotfill 8.1-p.43
    \par\noindent elementary embedding \dotfill 8.1-p.43
    \par\noindent elementary extension \dotfill 8.1-p.40
    \par\noindent elementary substructure \dotfill 8.1-p.40
    \par\noindent end extension \dotfill 11.12-p.63
    \par\noindent eventually different \dotfill 2.37-p.13
    \par\noindent exact cover \dotfill 1.33-p.8
    \par\noindent existential \dotfill 10.7-p.59
    \par\noindent extension \dotfill 8.1-p.40
    \par\noindent FIN \dotfill 11.10-p.63
    \par\noindent finitely axiomatizable \dotfill 7.4-p.36
    \par\noindent finitely generated \dotfill 8.30-p.47
    \par\noindent finitely realizable \dotfill 1.15-p.6
    \par\noindent finitely satisfiable \dotfill 7.1-p.31
    \par\noindent formulas \dotfill 7.1-p.25
    \par\noindent free \dotfill 7.1-p.26
    \par\noindent full arithmetic \dotfill 8.35-p.48
    \par\noindent graph \dotfill 1.30-p.8
    \par\noindent Hartog's ordinal \dotfill 6.23-p.23
    \par\noindent Henkin constant \dotfill 7.1-p.32
    \par\noindent Henkin sentence \dotfill 7.1-p.32
    \par\noindent Henkin \dotfill 7.1-p.32
    \par\noindent Herbrand \dotfill 10.8-p.59
    \par\noindent hereditarily finite \dotfill 11.1-p.62
    \par\noindent HF \dotfill 11.1-p.62
    \par\noindent incompleteness theorem \dotfill 11.1-p.62
    \par\noindent inconsistent \dotfill 10.1-p.57
    \par\noindent indecomposable \dotfill 5.23-p.20
    \par\noindent infix \dotfill 7.1-p.26
    \par\noindent interpretable \dotfill 11.18-p.64
    \par\noindent isomorphic \dotfill 8.1-p.42
    \par\noindent isomorphism \dotfill 3.4-p.15
    \par\noindent language of pure equality \dotfill 7.1-p.25
    \par\noindent lexicographical order \dotfill 6.4-p.22
    \par\noindent limit ordinals \dotfill 5.2-p.17
    \par\noindent linear order \dotfill 1.24-p.7
    \par\noindent logical complexity \dotfill 7.1-p.29
    \par\noindent logical symbols \dotfill 1.1-p.4
    \par\noindent logical symbols \dotfill 7.1-p.25
    \par\noindent logically equivalent \dotfill 1.1-p.4
  \par\noindent {\L}os-Vaught test \dotfill 8.7-p.45
    \par\noindent maximal \dotfill 1.21-p.6
    \par\noindent Maximality Principle \dotfill 4.1-p.16
    \par\noindent model \dotfill 7.1-p.27
    \par\noindent Modus Ponens \dotfill 10.1-p.56
    \par\noindent nonlogical symbols \dotfill 1.1-p.4
    \par\noindent nonlogical symbols \dotfill 7.1-p.25
    \par\noindent nonstandard model \dotfill 8.32-p.47
    \par\noindent nor \dotfill 1.8-p.5
    \par\noindent ordered pair \dotfill 2.3-p.11
    \par\noindent Ordinal arithmetic \dotfill 5.16-p.19
    \par\noindent ordinal \dotfill 5.1-p.17
    \par\noindent PA \dotfill 11.40-p.67
    \par\noindent partial order \dotfill 1.24-p.7
    \par\noindent partial recursive function \dotfill 9.1-p.50
    \par\noindent Peano arithmetic \dotfill 11.40-p.67
    \par\noindent postfix \dotfill 7.1-p.26
    \par\noindent prefix \dotfill 7.1-p.26
    \par\noindent prenex normal form \dotfill 10.4-p.58
    \par\noindent proof \dotfill 10.1-p.56
    \par\noindent proper classes \dotfill 2.1-p.11
    \par\noindent propositional letters \dotfill 1.1-p.4
    \par\noindent propositional sentences \dotfill 1.1-p.4
    \par\noindent realizable \dotfill 1.15-p.6
    \par\noindent reasonable proof predicate \dotfill 12.5-p.70
    \par\noindent recursive model \dotfill 10.25-p.61
    \par\noindent recursive set \dotfill 9.11-p.51
    \par\noindent recursive \dotfill 9.1-p.50
    \par\noindent recursively axiomatizable \dotfill 10.13-p.60
    \par\noindent recursively enumerable \dotfill 9.17-p.52
    \par\noindent regular \dotfill 6.10-p.22
    \par\noindent relatively prime \dotfill 11.24-p.65
    \par\noindent representable \dotfill 12.1-p.68
    \par\noindent restriction \dotfill 2.6-p.11
    \par\noindent second incompleteness \dotfill 12.1-p.68
    \par\noindent semantics \dotfill 1.1-p.4
    \par\noindent singular \dotfill 6.10-p.22
    \par\noindent Skolemization \dotfill 10.6-p.59
    \par\noindent soundness theorem \dotfill 10.1-p.57
    \par\noindent splits \dotfill 1.34-p.9
    \par\noindent state transition function \dotfill 9.1-p.50
    \par\noindent strict order \dotfill 1.24-p.7
    \par\noindent strongly undecidable \dotfill 11.17-p.64
    \par\noindent strongly undecidable \dotfill 11.19-p.64
    \par\noindent structure \dotfill 7.1-p.27
    \par\noindent substitution \dotfill 7.1-p.26
    \par\noindent substructure \dotfill 8.1-p.40
    \par\noindent successor ordinals \dotfill 5.2-p.17
    \par\noindent superstructure \dotfill 8.1-p.40
    \par\noindent syntax \dotfill 1.1-p.4
    \par\noindent terms \dotfill 7.1-p.25
    \par\noindent theory of groups \dotfill 7.1-p.25
    \par\noindent theory of partially ordered sets \dotfill 7.1-p.25
    \par\noindent Transfinite Recursion \dotfill 5.9-p.18
    \par\noindent transitive \dotfill 5.1-p.17
    \par\noindent transversal \dotfill 1.32-p.8
    \par\noindent truth evaluation \dotfill 1.1-p.4
    \par\noindent Tuckey's Lemma \dotfill 4.1-p.16
    \par\noindent Turing machine \dotfill 9.1-p.49
    \par\noindent twin prime conjecture \dotfill 8.32-p.47
    \par\noindent uncountable \dotfill 2.18-p.12
    \par\noindent universal \dotfill 10.6-p.59
    \par\noindent validity \dotfill 1.1-p.4
    \par\noindent Well-ordering Principle \dotfill 4.1-p.16
    \par\noindent wellorder \dotfill 3.1-p.15
    \par\noindent ZF \dotfill 2.1-p.10
    \par\noindent ZFC \dotfill 2.1-p.10
    \par\noindent ZFCFIN \dotfill 11.40-p.67
    \par\noindent Zorn's Lemma \dotfill 4.1-p.16
    \par\noindent \prf(x) \dotfill 12.4-p.69
 \immediate\closeout\dexno 
 \end{document}